\documentclass[11pt,a4paper,twoside]{article}
\title{\bf On real growth and run-off companies in insurance ruin theory}
\author{{\sc Harri Nyrhinen}}
\date{University of Helsinki, \today}

\usepackage{latexsym}
\usepackage{amsmath}
\usepackage{amssymb}
\usepackage{chicago}
\usepackage{amsbsy}

\newcommand{\halmos}{\vspace{3mm} \hfill \mbox{$\Box$}}

\newcommand{\reals}{{\mathbb R}}
\newcommand{\nat}{{\mathbb N}}

\newcommand{\rate}{{\mathfrak r}}

\newtheorem{Th}{Theorem}[section]
\newtheorem{Prop}{Proposition}[section]
\newtheorem{Lm}{Lemma}[section]

\newcommand{\Prob}{\mathbb{P}} 
\newcommand{\Exp}{\mathbb{E}} 
\newcommand{\Var}{{\rm Var}\,}

\begin{document}
\renewcommand{\baselinestretch}{1.0}
\maketitle

\renewcommand{\baselinestretch}{1.05}


\renewcommand{\baselinestretch}{1.05}

\begin{abstract}
We study solvency of insurers in a comprehensive model where various economic factors affect the capital developments of the companies. The main interest is
in the impact of real growth to ruin probabilities. The volume of the business is allowed to increase or decrease.
In the latter case, the study is focussed on run-off companies. Our main results give sharp asymptotic estimates for infinite time ruin probabilities.
\end{abstract}

%
%
%
%
\section{Introduction}\label{intro}
\setcounter{equation}{0}

Let $\{U_n: n = 0, 1, 2,\ldots\}$ be a real valued stochastic process
which describes the development of the capital of an insurance company.
Let the initial capital $U_0$ be a positive constant $u$.
The time of ruin $T = T_u$ is by definition,
\begin{equation}\label{s1f1}
T =
\begin{cases}
\inf\{n\in\nat\, ;\, U_n<0\}\cr
           \infty\mbox{ if }U_n\ge 0 \mbox{ for every } n\in\nat.
\end{cases}
\end{equation}
We are interested in the ruin probability $\Prob(T<\infty)$
for large $u$.

In classical risk theory, the capital development is
described by means of a random walk. The increments model yearly net incomes of the company, namely, differences between the premiums and the claims. Typically, the process $\{U_n\}$ has a linear drift to infinity.
In recent years, a lot of attention has been paid to models which allow 
economic factors to affect the capital development. Examples of such factors are inflation and returns on the investments. A key feature is that they cause multiplicative drifts to the capital process. It is nowadays understood
that the economic factors have a crucial impact
to ruin probabilities.

Real growth is a further economic factor which is motivated in
insurance context in applied studies of \shortciteN{TPJR82} and \shortciteN{CDTPMP94}. The feature is modelled as a trend in the numbers of claims in a multiplicative way. Periods of consecutive increase of the business volume may be very long. This phenomenon has been seen in the car insurance simply because the number of the cars has increased for a long time.
We will also study models in which the volume is drifting to zero. The main application in our mind is the case where solvency control is based on break-up basis.
Then it is assumed that the writing of new business is ceased so that the company is in the
run-off state. It has still to pay out compensations associated with the claims which have occurred but have not yet been settled, possibly not even notified. This liability of the company is a common feature in insurance contracts. The time to the final payment can be decades. In this context, ruin occurs if the capital and the investment income together do not suffice for the compensations.
These viewpoints are discussed in \shortciteN{TPHBMPJRMR89}, sections 1.2, 3.1.4 and 5.1. Detailed mathematical descriptions of the structure of the payment process in related models can be found in
\shortciteN{JR84}, \shortciteN{ICA88} and \shortciteN{RN93}.

Our purpose is  to give insight into risks associated with the cases where long drifts in the business volume are possible. To this end, we focus on the models which allow the volume to increase or decrease forever. This feature should be taken as an approximation of reality.

Real growth is not much studied in insurance ruin theory. We gave in \shortciteN{HN10} crude estimates for finite time ruin probabilities in this context. The focus was on increasing volumes. The results indicate that real growth
is then a substantial risk factor. In the present paper, the objective is to sharpen the view by deriving the asymptotic form for the ruin probability.
If the business volume is drifting to zero then new phenomena take place. It turns out that then ruin is likely to be caused by a single claim
at a late time point. This can be seen as a theoretical description of the run-off risk.

Asymptotic estimates for ruin probabilities in related earlier studies are largely based on the results of \shortciteN{CG91}.
The conclusion is that
\begin{equation}\label{s1f41}
\Prob(T < \infty) = (1+o(1)) C u^{-\kappa}, \quad u\to\infty,
\end{equation}
where $C$ and $\kappa$ are constants. A notable feature is that the key parameter $\kappa$ is merely determined by the economic factors. 
A survey of applications to ruin theory is given in
\shortciteN{JP08}.
Our model does not fit to this framework exactly, but
we will end up with estimates of $\Prob(T < \infty)$ by means of suitable approximations of the capital process.

The rest of the paper is organized as follows. The model for the capital development is described in Section \ref{Mr}. The main results are stated in Section \ref{sec20}. They are illustrated by means of
examples in sections \ref{Examples} and \ref{Sim}. The proofs are given in Section \ref{Proofs}.
%
%
%

\section{The model}\label{Mr}
\setcounter{equation}{0}

\vspace*{.05cm}
We describe in this section the basic structure of our model and give technical conditions which are assumed to be satisfied throughout
the paper. The model will be to large extent the same as in \shortciteN{HN10}.

We begin by describing the main variables and parameters of the model.

\subparagraph*{Numbers of claims}
Associated with year $n$, write
\begin{eqnarray}
N_n &=& \mbox{the accumulated number of claims up to year } n,
\nonumber\\
K_n &=& N_n -N_{n-1},
\mbox{ the number of claims in year }n,\nonumber\\
\lambda &=& \mbox{the basic level of the mean of the number of claims},\nonumber\\
\xi_n &=& \mbox{the mixing variable describing fluctuations
in the numbers of claims}.\nonumber
\end{eqnarray}
We assume that conditionally, 
given $\xi_1,\ldots,\xi_{n}$, the variables
$K_1,\ldots,K_n$ are independent and $K_k$
has the Poisson distribution with the parameter
$\lambda \xi_k$ for $k = 1, \ldots, n$.
The drift in the business volume will be modelled as a part of the mixing variables. Details will be given in subsequent sections.

\subparagraph*{Total claim amounts}
Let
\begin{eqnarray}
X_n &=& \mbox{the total claim amount in year }n,\nonumber\\
Z_j &=& \mbox{the size of the }j\mbox{th claim
in the inflation-free economy},\nonumber\\
m_Z &=& \mbox{the mean of the claim size in the inflation-free economy},\nonumber\\
i_n &=& \mbox{the rate of inflation in year }n.\nonumber
\end{eqnarray}
We consider the model where
\begin{eqnarray}\label{s2f0a}
X_n &=& (1+i_1)\cdots(1+i_{n})\sum_{j=N_{n-1}+1}^{N_n} Z_{j}.
\end{eqnarray}

\subparagraph*{Premiums}
For year $n$, write
\begin{eqnarray}
P_n &=& \mbox{the premium income}.\nonumber
\end{eqnarray}
The structure of $P_n$ will be specified in subsequent sections.

\subparagraph*{The transition rule}
We next describe the development of the capital in time. Let
\begin{eqnarray}
U_n &=& \mbox{the capital at the end of year }n,\nonumber\\
r_n &=& \mbox{the rate of return on the investments
in year }n.\nonumber
\end{eqnarray}
Let $U_0 = u >0$ be the deterministic initial capital
of the company. We define
\begin{eqnarray}\label{s2f2}
U_n &=& (1+r_n)(U_{n-1}+P_n-X_n).
\end{eqnarray}
This transition rule corresponds to the case where the premiums and the claims are all paid at the beginning of the year.
It would also be natural to define
\begin{eqnarray}\label{s2f2late}
U_n = (1+r_n)(U_{n-1}+P_n)-X_n.
\end{eqnarray}
Then the premiums are paid at the beginning and the claims at the end of the year. The reality is probably somewhere between
\eqref{s2f2} and \eqref{s2f2late}. We will assume \eqref{s2f2} in the sequel but it should not be very different to analyse \eqref{s2f2late}.

\subparagraph*{Technical specifications and assumptions}

We end the description by specifying the dependence structure
and other technical features of the model.
All the random variables below are assumed to be
defined on a fixed probability space $(\Omega, {\cal F}, \Prob)$.

We begin by giving a detailed mathematical description for the total claim amounts
in the inflation-free economy. For year $n$, denote this quantity by $V_n$, that is,
\begin{equation}\label{s2f20}
V_n = \sum_{j=N_{n-1}+1}^{N_n} Z_j.
\end{equation}
The distributions of the $N$- and $K$-variables depend on $\xi$-variables. We assume that
$\lambda > 0$ and $\Prob(\xi_n >0) = 1$
for every $n\in\nat$. Denote by $F_n$ the joint distribution function of the vector
$(\xi_1,\ldots,\xi_n)$. We assume that
for every $h_1,\ldots,h_n\in \nat\cup\{0\}$
and for every Borel set $C\subseteq\reals^{n}$,
\begin{eqnarray}\label{s2f0b1}
&&\Prob\left(K_1 = h_1,\ldots,K_n = h_n,
(\xi_1,\ldots,\xi_n)\in C\right)\\
&=& \int_{(y_1,\ldots,y_n)\in  C}
\prod_{k=1}^n e^{-\lambda y_k}
\frac{(\lambda y_k)^{h_k}}{h_k!}\,
dF_n(y_1,\ldots,y_n).\nonumber
\end{eqnarray}
The claim sizes
$Z, Z_1, Z_2,\ldots$
are assumed to be i.i.d.
We also assume that they are
independent of the numbers of claims in all respects. Let $F_Z$ be the distribution function of $Z$, and let $F_Z^{h*}$ be the $h$th convolution power of $F_Z$.
We assume that
for every $h_1,\ldots,h_n\in \nat\cup\{0\}$
and $v_1,\ldots,v_n\in\reals$,
and for every Borel set $C\subseteq\reals^{n}$,
\begin{eqnarray}\label{s2f0c}
&&\Prob\left(V_1\le v_1,\ldots,V_n\le v_n,
K_1 = h_1,\ldots,K_n = h_n,
(\xi_1,\ldots,\xi_n)\in C\right)\\
&=& \Prob\left(K_1 = h_1,\ldots,K_n = h_n,
(\xi_1,\ldots,\xi_n)\in C\right)
\prod_{k=1}^n F_Z^{h_k*}(v_k).\nonumber
\end{eqnarray}
We refer to \shortciteN{JG97} for more information about
mixed Poisson distributions.

Consider now the other parts of the model.
Concerning inflation and the returns on the investments, we take
$(i, r), (i_1,r_1), (i_2,r_2),\ldots$
to be an i.i.d. sequence of random vectors, and these vectors are assumed
to be independent of $\xi$-, $K$- and $Z$-variables.
For the supports of $Z$, $i$ and $r$, we assume that
$\Prob(Z > 0) = 1,\Prob(i>-1) = 1$ and
$\Prob(r>-1) = 1$.

\section{Main results}\label{sec20}

\setcounter{equation}{0}
Let the model be as described in Section \ref{Mr} and let the time of ruin $T$  be as in \eqref{s1f1}.
It is convenient to consider a discounted version of the process $\{U_n\}$. Write
\begin{equation}\label{late1}
A = \frac{1+i}{1+r}\quad\mbox{and}\quad A_n = \frac{1+i_n}{1+r_n}
\end{equation}
for $n\in\nat$, and let
\begin{equation}\label{pf10b}
B_n = V_n-\frac{P_n}{(1+i_1)\cdots(1+i_n)}
\end{equation}
where $V_n$ is as in \eqref{s2f20}.
Write further
\begin{equation}\label{pf10b1}
Y_n = \sum_{k=1}^n A_1\cdots A_{k-1} (1+i_k) B_k.
\end{equation}
By dividing $U_n$ by $(1+r_1)\cdots(1+r_n)$, it is seen that
the time of ruin can be expressed as
\begin{equation}\label{pt3}
T =
\begin{cases}
\inf\{n\in\nat\, ;\, Y_n>u\}\cr
           \infty\mbox{ if }Y_n\le u \mbox{ for every } n\in\nat.
\end{cases}
\end{equation}
The ruin probability can also be defined by means of
\begin{equation}\label{p14f10b1}
\bar{Y} := \sup\{Y_n; n = 1, 2,\ldots\}.
\end{equation}
Namely, $\Prob(T<\infty) = \Prob(\bar{Y} > u)$.

\subsection{Background results}

We present in this section mathematical tools which are of general interest and which are needed in our study.
The symbol $=_L$ will mean equality of probability laws and the notation $a^+$ the positive part of $a\in\reals$.
\begin{Th}\label{basic}
Let $(M,Q)$ be a two-dimensional random vector. Assume that $\Prob(M\ge 0) = 1$
and that for some $0<\kappa < \alpha$,
$$\Exp(M^{\kappa}) = 1\quad\mbox{and}\quad
\Exp(M^{\alpha}),\, \Exp(|Q|^{\alpha}) < \infty.$$
Assume further that the conditional law of $\log M$, given $M\not= 0$, is non-arithmetic. Then there exists a random variable $R$ which satisfies the random equation
\begin{equation}\label{s14m3}
R =_L Q + M\max(0,R),\quad R \mbox{ independent of } (M,Q).
\end{equation}
Furthermore, if $R$ satisfies \eqref{s14m3} then
\begin{equation}\label{s14m3a}
\lim_{u\to\infty}
u^{\kappa} \Prob(R>u) = C
\end{equation}
where
\begin{equation}\label{s14m3b}
C = \frac{\Exp\left(((Q + M\max(0,R))^+)^{\kappa} - ((MR)^+)^{\kappa}\right)}
{\kappa m}
\end{equation}
and $m = \Exp(M^{\kappa}\log M)$.
\end{Th}
The proof of the result can be found in \shortciteN{CG91}, Theorem 6.2. Theorem \ref{basic} has been applied directly to ruin theory, for example, in \shortciteN{HN01}. In the present model,
we need also the following approximation scheme.

\begin{Lm} \label{cvLm}
Let $\{{\cal Y}_n\}$, $\{{\cal Y}_{n1}\}$ and $\{{\cal Y}_{n2}\}$ be stochastic processes
such that
\begin{equation}\label{s14m9a}
{\cal Y}_n = {\cal Y}_{n1} + {\cal Y}_{n2},\quad n = 1, 2,\ldots.
\end{equation}
Write
\begin{equation}\label{s14m10}
\bar{{\cal Y}} = \sup\{{\cal Y}_n; n\in\nat\}\quad
\mbox{and}\quad
\bar{{\cal Y}}_1 = \sup\{{\cal Y}_{n1}; n\in\nat\},
\end{equation}
and let ${\cal Y}_{02} \equiv 0$. Assume that there exists
$\kappa\in (0,\infty)$, $\alpha\in (\kappa,\infty)$ and
$\delta\in (0, 1-\frac{\kappa}{\alpha})$ such that
\begin{equation}\label{s14m10aa}
\Exp\left(\left|{\cal Y}_{n2} -
{\cal Y}_{n-1,2}\right|^{\alpha}\right) < \infty,\quad n = 1, 2,\ldots,
\end{equation}
\begin{equation}\label{s14m10a}
\limsup_{n\to\infty} n^{-1}\log
\Exp\left(\left|{\cal Y}_{n2} -
{\cal Y}_{n-1,2}\right|^{\alpha}\right) < 0,
\end{equation}
\begin{equation}\label{s14m3a0}
\liminf_{u\to\infty}
(\log u)^{-1} \log \Prob(\bar{{\cal Y}}_1>u) \ge -\kappa
\end{equation}
and
\begin{equation}\label{s14m3a1}
\Prob\left(\bar{{\cal Y}}_1>u(1+u^{-\delta})\right) = (1+o(1)) \Prob(\bar{{\cal Y}}_1>u),
\quad u\to\infty.
\end{equation}
Then
\begin{equation}\label{s14m3a2}
\Prob(\bar{{\cal Y}}>u) =  (1+o(1)) \Prob(\bar{{\cal Y}}_1>u),\quad u\to\infty.
\end{equation}
\end{Lm}
A possible way to apply the lemma in our model is to take
$\{{\cal Y}_n\} = \{Y_n\}$ and to find out a suitable process
$\{{\cal Y}_{n1}\}$
such that Theorem \ref{basic} can be used to conclude that
\begin{equation}\label{s14mb10}
\lim_{u\to\infty}
u^{\kappa}
\Prob(\bar{{\cal Y}}_1>u) = C
\end{equation}
for some $\kappa >0$ and $C>0$. Then \eqref{s14m3a0} and
\eqref{s14m3a1} are automatically satisfied.
If also \eqref{s14m10aa} and \eqref{s14m10a} hold then we obtain an estimate for the ruin probability.

The last result gives descriptions of tails
associated with compound distributions.
Let $\eta, \eta_1, \eta_2,\ldots$ be i.i.d
random variables, and let ${\cal V}_0 \equiv 0$ and
${\cal V}_n = \eta_1+\cdots +\eta_n$
for $n\in\nat$. Denote by $\Lambda_{\eta}$ the cumulant generating function of $\eta$, that is,
$\Lambda_{\eta}(\alpha) = \log \Exp\left(e^{\alpha \eta}\right)$
for $\alpha\in\reals$.
Let ${\cal W}$ and ${\cal N}$ be independent random variables which are also independent of $\eta$-variables.
Assume that $\Prob({\cal N}\in \nat\cup\{0\}) = 1$.
Recall that a function $f:(0,\infty)\to (0,\infty)$ is regularly varying if there exists $\gamma\in\reals$
such that for every $x > 0$,
\begin{equation}\label{rwlm21a}
\lim_{t\to\infty} \frac{f(tx)}{f(t)} = x^{\gamma}.
\end{equation}

\begin{Lm}\label{rwlm2}
Assume that the distribution
of $\eta$ is non-arithmetic and that
\begin{equation}\label{rwlm21}
\lim_{n\to\infty} n^{-1} \log \Prob({\cal N} = n) = -\upsilon
\end{equation}
where $\upsilon\in (0,\infty)$.
Assume further that there exists $\rate\in (0,\infty)$ such that
$\Lambda_{\eta}(\rate) = \upsilon$, and that
$\Lambda_{\eta}(\alpha)$ and $\Exp(e^{\alpha {\cal W}})$
are finite for some $\alpha >\rate$. Then
\begin{equation}\label{rwlm22}
\lim_{u\to\infty} u^{-1} \log
\Prob({\cal V}_{{\cal N}}+{\cal W} > u) = -\rate.
\end{equation}
Further, $\Lambda'_{\eta}(\rate) > 0$, and if
$\mu =  1/\Lambda'_{\eta}(\rate)$,
then for every $\varepsilon > 0$, there exists
$\varepsilon' > 0$ such that
\begin{equation}\label{rwlm23}
\Prob\left({\cal N}/u\in
[\mu-\varepsilon,\mu+\varepsilon]\,\left|\,
\right.
{\cal V}_{\cal N}+{\cal W} > u\right)
= 1+O(e^{-\varepsilon'u}),
\quad u\to\infty.
\end{equation}
If in addition,
\begin{equation}\label{rwlm25}
\Prob({\cal N} = n) = (1+o(1)) f(n) e^{-n\upsilon},\quad n\to\infty,
\end{equation}
where $f$ is regularly varying then
\begin{equation}\label{rwlm27}
\Prob({\cal V}_{{\cal N}}+{\cal W} > u)
= (1+o(1))\,
\frac{\Exp(e^{\rate {\cal W}}) \mu}{\rate}
f\left(\mu u\right) e^{-\rate u},\quad u\to\infty.
\end{equation}
\end{Lm}
The last estimate \eqref{rwlm27} is closely related to
\shortciteN{PEMMJT85} and
\shortciteN{JT85}. The main difference is that we allow negative values for $\eta$-variables.

\subsection{The case of increasing volumes}\label{sec21}

We give in this section estimates for ruin probabilities in the case where the business volume has a tendency to increase. We begin with some further specifications of our model. Recall the descriptions of the numbers of claims
$K_n$  from Section \ref{Mr}. Associated with year $n$, write
\begin{eqnarray}
g_n &=& \mbox{the rate of real growth},\nonumber\\
q_n &=& \mbox{the structure variable describing
short term fluctuations in the numbers of claims}.\nonumber
\end{eqnarray}
We take
\begin{equation}\label{s2gnot00}
\xi_n = (1+g_1)\cdots (1+g_n) q_n,\quad n\in\nat.
\end{equation}
Assume that
$(g,q), (g_1,q_1), (g_2,q_2),\ldots$
is an i.i.d. sequence of random vectors which is independent of the other variables of the model.
Assume also that $g$ and $q$ are independent and that
\begin{equation}\label{s2gnot0}
\Prob(g=0) < 1,\quad\Prob(g>-1) = 1,\quad
\Prob(q > 0) = 1\quad\mbox{ and }
\quad\Exp(q) = 1.
\end{equation}
The premium $P_n$ is supposed to have the form
\begin{equation}\label{s2f1}
P_n = (1+s) \lambda m_Z (1+g_1)\cdots(1+g_{n})
(1+i_1)\cdots(1+i_{n})
\end{equation}
where $s>0$ is the safety loading coefficient. Then
$Y_n$ of \eqref{pf10b1} has the form
\begin{equation}\label{pf10b1copy}
Y_n = \sum_{k=1}^n A_1\cdots A_{k-1} (1+i_k)
[V_k - (1+s) \lambda m_Z (1+g_1)\cdots(1+g_{k})].
\end{equation}
The above structure is suggested in \shortciteN{CDTPMP94}.

Define the functions $\Lambda_A, \Lambda_g$, $\Lambda_1$:
$\reals\to\reals\cup\{\infty\}$ by
\begin{eqnarray}
\Lambda_A(\alpha) &=& \log\Exp\left(A^{\alpha}\right),\label{s2f3}\\
\Lambda_g(\alpha) &=& \log\Exp\left((1+g)^{\alpha}\right),\label{s2f3a}\\
\Lambda_1(\alpha) &=& \Lambda_A(\alpha) + \Lambda_g(\alpha).\label{s2f3b}
\end{eqnarray}
Observe that $\Lambda_A, \Lambda_g$ and $\Lambda_1$ are cumulant generating functions so that they are convex.
Define the parameters $\rate_1$ and $\beta_1$ by
\begin{eqnarray}
\rate_1 &=& \sup\{\alpha\ge 0; \Lambda_1(\alpha) \le 0\}\label{s14m5}.
\end{eqnarray}
and
\begin{eqnarray}\label{pf1422a}
\beta_1 = \sup\{\alpha\in\reals\, |\,
\Lambda_1(\alpha),\, \Exp((1+i)^{\alpha}),\,\Exp(Z^{\alpha}),\,
\Exp(q^{\alpha}) < \infty\}\in [0,\infty].
\end{eqnarray}
Write
\begin{equation}\label{s2f3bb}
D = A(1+g)\quad\mbox{and}\quad D_n = A_n(1+g_n).
\end{equation}

\begin{Th} \label{mainthm1}
Let the model be as described above.
Assume that $\Exp(\log(1+g))\ge 0$ and that
$\beta_1 > 0$ and $\rate_1\in (0,\beta_1)$.
Assume further that $\Exp(Z^{\alpha})<\infty
\mbox{ for some } \alpha > 1$ and that the distribution of $\log D$ has a non-trivial absolutely continuous component. Let
\begin{equation}\label{s2f15B}
Q = (1+i)(1+g)\lambda m_Z(q-(1+s))
\quad\mbox{and}\quad M = D.
\end{equation}
Then $(Q,M)$ satisfies the conditions of Theorem \ref{basic}
with $\kappa = \rate_1$. If $R$ satisfies
\eqref{s14m3} and if $C$ is the constant of \eqref{s14m3b} then
\begin{eqnarray}\label{s2f15}
&&\lim_{u\to\infty}u^{\rate_1} \Prob(T < \infty) = C.
\end{eqnarray}
Furthermore, $C$ is strictly
positive if and only if $\Prob(q>1+s) > 0$.
\end{Th}
We apply Lemma \ref{cvLm} in the proof of Theorem \ref{mainthm1} by taking ${\cal Y}_n = Y_n$ and
\begin{equation}\label{pf1419}
{\cal Y}_{n1} = \sum_{k=1}^n
D_1\cdots D_{k-1} (1+i_k)(1+g_k) \lambda m_Z(q_k-(1+s)).
\end{equation}
Indeed, we show that
$\Prob(T<\infty) = (1+o(1))
\Prob(\bar{{\cal Y}}_{1}>u)$ as $u\to\infty$
where $\bar{{\cal Y}}_{1}$ is as in \eqref{s14m10}.
Real growth and inflation appear equally in
\eqref{pf1419} so that their impacts to ruin probabilities
are similar. However, a feature caused by non-degenerate real growth is that
the claim sizes $Z_1, Z_2,\ldots$
only contribute to ${\cal Y}_{n1}$ and to estimate
\eqref{s2f15}
via the mean $m_Z$. The claim numbers have
a more drastic effect via the structure variable $q$.
Even the positivity of $C$ in \eqref{s2f15}
depends on the support of $q$.

\subsection{The case of decreasing volumes}\label{sec21b}

We consider in this section ruin probabilities of run-off companies by allowing the business volume to drift to zero. The structure of the model will be as in Section \ref{Mr} but we now drop the premiums from
the considerations by taking
$P_n=0$ for each $n$.
The interpretation is that no new insurance contracts are made after year 0.
We assume that the company has liabilities from the past
so that it has to pay compensations associated with the claims which have occurred but have not yet been settled.
In this context, it is natural to interpret $K_n$ as the number of payments or as the number of reported claims in year $n$. We discuss the latter case in detail in Section \ref{Examples}.

We dropped the premiums so that the variable $Y_n$ of \eqref{pf10b1} is non-negative and
\begin{eqnarray}\label{s2f3p}
\bar{Y} = \sum_{n=1}^{\infty}
A_1\cdots A_{n-1} (1+i_n) V_n\quad \mbox{ a.s}.
\end{eqnarray}
The model for the $\xi$-variables of Section \ref{sec21}
is too simple from the applied point of view.
To get a suitable generalization, define the function
$\Lambda_{\xi}: \reals\to\reals\cup\{\pm\infty\}$ by
\begin{eqnarray}\label{s2f3l}
\Lambda_{\xi}(\alpha) &=&
\limsup_{n\to\infty} n^{-1}
\log\Exp\left(\xi_n^{\alpha}\right).
\end{eqnarray}
Then $\Lambda_{\xi}$ is convex.
We will work under the following hypotheses $(H1)-(H2)$. 

\begin{tabbing}
\=HHHHHHH
\=Assumption\kill
\>$(H1)$\>$\lim_{\alpha\to 0+}
\Lambda_{\xi}(\alpha) = 0\quad\mbox{ and }\quad
\lim_{\alpha\to \infty}
\Lambda_{\xi}(\alpha) = -\infty.$
\end{tabbing}
\begin{tabbing}
\=HHHHHHH
\=Assumption\kill
\>$(H2)$\>$\mbox{For }\alpha = 1,\,\,
\eqref{s2f3l}
\mbox{ holds as the limit}.$
\end{tabbing}

Basic facts concerning the numbers of claims are given in the following result. Sharper but more technical descriptions are stated in Lemma \ref{cnLm} in Section
\ref{Proofs}.

\begin{Prop}\label{Kprop}
Let the model be as described above.
Assume that $(H1)-(H2)$ are satisfied. Then
as $n\to\infty$,
\begin{eqnarray}
\Prob(K_n=1) &=& (1+o(1)) \lambda\Exp(\xi_n),\label{cnLm09}\\
\Prob(K_n\ge 2) &=& o(1) \Prob(K_n=1),\label{cnLm09x}
\end{eqnarray}
and
\begin{eqnarray}
\lim_{n\to\infty} n^{-1}
\log \Prob(K_n=1) &=& \Lambda_{\xi}(1),\label{cnLm00}\\
\lim_{n\to\infty} n^{-1}
\log \Prob(\xi_n\ge 1) &=& -\infty.
\label{cnLm11}
\end{eqnarray}
\end{Prop}

The descriptions of the proposition illustrate hypotheses
$(H1)-(H2)$. First observe that
 $\Lambda_{\xi}(1)\in (-\infty,0)$ by $(H1)$. By
\eqref{cnLm09x} and \eqref{cnLm00}, the probability
$\Prob(K_n\ge 1)$ tends to zero. This is a natural requirement for run-off companies. Limit \eqref{cnLm11} means that the random Poisson parameter
$\lambda\xi_n$ has a strong tendency to be below the basic level $\lambda$. Limit \eqref{cnLm00} also shows that $\Prob(K_n\ge 1)$ is positive for every $n$. This feature has a more theoretical nature but it may be viewed as an approximation of the reality in the case of long
delays in claims settlements. Theorem \ref{mainthm2} below indicates that ruin is likely to occur rather quickly so that the feature is perhaps not that critical from the applied point of view.

Consider now ruin probabilities.
Define the function $\Lambda_2:\reals\to\reals\cup\{\infty\}$ by
\begin{eqnarray}\label{s2f3c}
\Lambda_2(\alpha) &=&
\Lambda_A(\alpha) + \Lambda_{\xi}(1).
\end{eqnarray}
Then $\Lambda_2$ is convex.
Define the parameters $\rate_2$ and $\beta_2$ by
\begin{eqnarray}
\rate_2 &=& \sup\{\alpha\ge 0; \Lambda_2(\alpha) \le 0\}\label{s14m6}
\end{eqnarray}
and
\begin{eqnarray}\label{pf1422aa}
\beta_2 = \sup\{\alpha\in\reals\, |\,
\Lambda_2(\alpha),\, \Exp((1+i)^{\alpha}),\,\Exp(Z^{\alpha}) < \infty\}\in [0,\infty].
\end{eqnarray}
We will assume below that $\rate_2\in (1,\beta_2)$.
Then $\Lambda_2(\rate_2) = 0$ so that under $(H1)$,
$$\Lambda_A(\rate_2) > 0\quad
\mbox{and}\quad\Lambda'_A(\rate_2) > 0.$$
Write
$\mu_2 = 1/\Lambda'_A(\rate_2)$.

We will apply Lemma \ref{basic} by taking
${\cal Y}_n = Y_n$ and
\begin{equation}\label{c14203a1}
{\cal Y}_{n1} = \sum_{k=1}^{n}
A_1\cdots A_{k-1} (1+i_k) V_k
{\bf 1}(K_k = 1, K_j= 0,\, \forall j\ge k+1).
\end{equation}
Then
\begin{equation}\label{c14203a}
\bar{{\cal Y}}_{1} = \sum_{n=1}^{\infty}
A_1\cdots A_{n-1} (1+i_n) V_n
{\bf 1}(K_n = 1, K_j= 0,\, \forall j\ge n+1)\quad
\mbox{a.s.}
\end{equation}

\begin{Th} \label{mainthm2}
Assume $(H1)-(H2)$ and that $\beta_2 > 1$
and $\rate_2\in (1,\beta_2)$.
Assume further that the distribution
of $\log A$ is non-lattice. Let $\bar{{\cal Y}}_{1}$ be as in \eqref{c14203a}. Then
\begin{eqnarray}
\lim_{u\to\infty}
(\log u)^{-1} \log \Prob(T<\infty) &=& -\rate_2\label{s142030},
\end{eqnarray}
\begin{eqnarray}
&&\Prob(T<\infty) = (1+o(1))\,
\Prob(\bar{{\cal Y}}_{1} > u)
\label{c14203at}\\
&=& (1+o(1))\,\sum_{n=1}^{\infty}
\Prob(A_1\cdots A_{n-1} (1+i) Z> u)\,
\Prob(K_n = 1),\quad u\to\infty,\label{c14203aa}
\end{eqnarray}
and for every $\varepsilon > 0$,
\begin{equation}\label{cnLm09latetyp}
\lim_{u\to\infty}
\Prob\left(T/\log u \in [\mu_2-\varepsilon,\mu_2+\varepsilon]
\,|\, T<\infty\right) = 1.
\end{equation}
If in addition,
\begin{equation}\label{cnLm09late}
\Prob(K_n = 1) = (1+o(1)) \lambda f(n)
e^{n\Lambda_{\xi}(1)},\quad n\to\infty,
\end{equation}
where $f$ is regularly varying then as $u\to\infty$,
\begin{eqnarray}
\Prob(T<\infty) &=& (1+o(1))\,
\frac{\Exp((1+i)^{\rate_2})
\Exp(Z^{\rate_2})e^{\Lambda_{\xi}(1)}\mu_2}{\rate_2}
\lambda f\left(\mu_2\log u\right) u^{Ì-\rate_2}.\label{c14203aal}
\end{eqnarray}
\end{Th}

Estimate \eqref{c14203at} shows that the tail
probabilities of $\bar{Y}$ of \eqref{s2f3p} and $\bar{{\cal Y}}_{1}$ of \eqref{c14203a} are asymptotically equivalent. This is suprising since the two variables are the same except that
$\bar{{\cal Y}}_{1}$ disregards a big part of the claims.
A possible intuitive interpretation is that
in order to get ruined, the company first looses a major part of its capital mainly because of bad returns on the investments, and the rest of the capital is lost by the very last claim. This phenomenon is somewhat strange but a dominance of a single claim has also been found in connection with heavy tailed claim sizes. We refer the reader to
\shortciteN{SACK96}. We do not assume heavy tails but it is worth to observe that late claims may be large because of high inflation. It could also be possible to find out different views by making use of alternative limiting procedures. For example, if we would allow $\lambda$ to increase with $u$ then also earlier claims could contribute meaningfully the ruin probability.

We can expect that the accuracy
of estimate \eqref{c14203at}
is not very good for moderate initial capitals.
To be accurate, the probabilities
$$\Prob(K_n=1)\quad\mbox{and}\quad \Prob(K_n\ge1)$$
should be close to each other at least for $n$ close to $\mu_2\log u$.
A large $\lambda$ easily violates this relation. We consider the problem quantitatively in Section \ref{Sim}.

Estimate \eqref{c14203aa} is connected with tails of compound distributions. To see this, write
\begin{eqnarray}\label{ss0}
p_n = \Prob(K_n = 1)\quad \mbox{and}\and\quad
p = \sum_{n=1}^{\infty} p_n.
\end{eqnarray}
It follows from \eqref{cnLm00} that $p\in (0,\infty)$.
Let $\rho$ be a random variable such that
\begin{equation}\label{ss00}
\Prob(\rho = n-1) = p_n/p,
\quad n\in\nat,
\end{equation}
and assume that $\rho$ is independent of everything else.
Write further $S_0 = 0$ and
\begin{equation}\label{ss01}
S_n = \log A_1 + \cdots +\log A_n,\quad n \in\nat.
\end{equation}
By \eqref{c14203at} and \eqref{c14203aa},
\begin{eqnarray}
\Prob(T < \infty) &=&
(1+o(1)) p\,\Prob(A_1\cdots A_{\rho} (1+i)Z > u)
\label{lateT}\\
&=&(1+o(1)) p\,\Prob(S_{\rho} + \log((1+i)Z) > \log u).
\nonumber
\end{eqnarray}
The last probability can be approximated by means of Lemma \ref{rwlm2}.

\section{An applied example}\label{Examples}
\setcounter{equation}{0}

Consider a run-off company which
has operated in the market in years
$-d,\ldots,-1, 0$
where $d\in\nat$. Let $u>0$ be the capital of the company at the end of year 0. The company does not make insurance contracts in the future. Thus ruin means that the capital and the investment income together do not suffice for the compensations to be paid in years
$1, 2,\ldots$.

A suitable way to model the future events is to associate with each claim the {\it reporting time}, that is,
the time at which the company receives the first information about the claim. The {\it reporting delay} is
the difference between the reporting and the occurrence time of the claim. We assume that the compensations are paid at the reporting times.

We will take the model of Section \ref{sec21} as the description of occurrences of claims in years $-d,\ldots,0$. There is no need to describe premiums or returns on the investments for the past years since their affects are accumulated into the initial capital $u$.
With year $m\in \{-d,\ldots, 0\},$
associate the structure variable
$q_m$, and assume that
$q, q_{-d},\ldots,q_{0}$
are i.i.d. random variables. Assume also that $\Prob(q > 0) = 1$ and $\Exp(q) = 1$.
Let $\pi_0 = 1$ and $\pi_{-d},\ldots,\pi_{-1}$ be positive constants which describe the observed levels of the business volume in the past years. We assume that in year $m$, claims have occurred according to a mixed Poisson process such that conditionally, given $q_m$, the intensity of the process has been
$\lambda \pi_m q_m$. For different years, the occurrence processes are assumed to be independent. The reporting delays are assumed to be i.i.d. random variables with the common distribution function $G$ with $G(0)=0$. Assume also that they are independent of everything else.
Inflation is assumed to affect such that the size of any reported claim in year $n\ge 1$ has the same distribution as
$$(1+i_1)\cdots(1+i_n) Z.$$

Fix $m\in \{-d,\ldots,0\}$ and consider claims which have occurred in year $m$.
The number of reported claims in year
$n\ge 1$ has a mixed Poisson distribution. The random Poisson parameter is $\lambda \pi_m b_{n-m} q_m$ where
\begin{equation}
b_{k} = \int_{0}^{1} (G(k+1-s)-G(k-s)) ds,
\quad k\in\nat.\nonumber
\end{equation}
A further useful fact is that conditionally,
given $q_{m}$, the numbers of reported claims in different years are independent.
We refer the reader to \shortciteN{JR84}, Section 2.3.1.

By the above discussion, the number of reported claims in year $n$ has a mixed Poisson distribution. The Poisson parameter is $\lambda\xi_n$ where
\begin{equation}
\xi_n = \sum_{m=-d}^0 \pi_m b_{n-m} q_m.\nonumber
\end{equation}
Our basic assumption \eqref{s2f0b1} is also satisfied.

Assume that $\Exp(q^{\alpha})<\infty$
for every $\alpha > 0$ and that
\begin{equation}\label{tailg}
1-G(x) = (1+o(1)) h(x) e^{-x\varphi},\quad x\to\infty,
\end{equation}
where $h$ is regularly varying
and $\varphi\in (0,\infty)$
is a constant. Requirement \eqref{tailg} is satisfied, for example, by every gamma distribution. It is easy to see that
$\Lambda_{\xi}(\alpha) = -\alpha\varphi$ for every $\alpha > 0$
and that also $(H2)$ is satisfied. Furthermore,
\eqref{cnLm09late} holds since
\begin{equation}\label{exf1}
\Prob(K_n=1) = (1+o(1)) \lambda
\frac{(e^{\varphi}-1)(1-e^{-\varphi})
\sum_{m=-d}^0 \pi_m e^{m\varphi}}{\varphi} h(n) e^{-n\varphi}.
\end{equation}

\section{A simulation example}\label{Sim}
\setcounter{equation}{0}

The asymptotic estimate of Theorem \ref{mainthm2}
disregards a lot of claims so that it is interesting to study its accuracy for moderate initial capitals $u$.
We accomplish this by means of simulation.
We also suggest an ad hoc method to estimate
efficiently ruin probabilities.

We begin by fixing the model to be considered.
Concerning the returns on the investments, we assume that
$\log (1+r)$ has a normal distribution. Denote by $m_r$ and $\sigma_r$
the mean and the standard deviation, respectively. The rate of inflation is a constant. Write in short $m_i = \log (1+i)$. The mixing variables will also be deterministic. We take
$\xi_n = e^{-n\varphi}$
where $\varphi\in (0,\infty)$ is a constant.
Finally, the claim size $Z$ will be exponentially distributed.

By the above specifications, we have
for $\alpha\in\reals$ and $n\in\nat$,
\begin{eqnarray}
\Lambda_A(\alpha) &=& (m_i-m_r)\alpha +
\sigma_r^2\alpha^2/2,\nonumber\\
\Lambda_2(\alpha) &=& \Lambda_A(\alpha) - \varphi,\nonumber\\
\Lambda_{\xi}(\alpha) &=& -\varphi\alpha,\nonumber\\
\Prob(K_n=1) &=& (1+o(1)) \lambda e^{-n\varphi}.\nonumber
\end{eqnarray}
Thus $f\equiv 1$ in Theorem \ref{mainthm2}.
The numeric values of the parameters will be
$$m_r = 0.1,\quad\sigma_r^2 = 0.1,\quad m_i = 0.05,
\quad \varphi = 0.1\quad
\mbox{and}\quad \Exp(Z) =1.$$
Then $\rate_2=2$ and the estimate of Theorem \ref{mainthm2} is
\begin{eqnarray}\label{rlate}
\Prob(T<\infty) = (1+o(1))
\frac{20}{3}\lambda u^{-2}.
\end{eqnarray}
In the following tables 5.1 and 5.2, we use notations
\begin{eqnarray}
\hat{E_1} &=&  \mbox{estimate }
\eqref{rlate} \mbox{ of the ruin probability
with }o(1)\mbox{ replaced by zero},\nonumber\\
\hat{E_2} &=& \mbox{the estimate of the ruin probability from simulation}.\nonumber
\end{eqnarray}
We had approximately 10 millions replications in simulation of each probability so that estimates $\hat{E}_2$ should be rather close to the true values. The quotient $\hat{E}_2/\hat{E}_1$ describes the accuracy of $\hat{E}_1$. In Table 5.1, it is rather close to one as it should be. In Table 5.2, $\lambda$ is larger
and the resulting $\hat{E}_2/\hat{E}_1$ is large so that the estimate of Theorem \ref{mainthm2} is inaccurate.

\hskip 0.5 pt

\leftline{\:\:\:\:\:\:\:\:\:\:\:\:\:\:\:\:\:\:\:\:\:\:\:\:\:\:\:\:\:\:Table 5.1, $\lambda = 0.1$}
\hskip 0.5 pt

\begin{tabular}{|c|c|c|c|c|c|}
\hline

$u$ & \footnotesize$\hat{E}_1$ & \footnotesize$\hat{E}_2$ &
\footnotesize$\hat{E}_2/\hat{E}_1$\\
\hline
\small$10$ & 6.7$\times 10^{-3}$ & $9.3\times 10^{-3}$ & 1.40\\
\hline
\small$50$ & 2.7$\times 10^{-4}$ & $3.2\times 10^{-4}$ & 1.20\\
\hline
\small$200$ & 1.7$\times 10^{-5}$ & $1.9\times 10^{-5}$ & 1.12\\
\hline
\end{tabular}

\hskip 0.5 pt

\bigskip
\leftline{\:\:\:\:\:\:\:\:\:\:\:\:\:\:\:\:\:\:\:\:\:\:\:\:\:\:\:\:\:\:Table 5.2, $\lambda = 100$}
\hskip 0.5 pt
\begin{tabular}{|c|c|c|r|c|c|}
\hline

$u$ & \footnotesize$\hat{E}_1$ & \footnotesize$\hat{E}_2$ &
\footnotesize$\hat{E}_2/\hat{E}_1$ & \footnotesize$\hat{E}_3$ &
\footnotesize$\hat{E}_3/\hat{E}_2$\\
\hline
\small$5\,000$ & 2.7$\times 10^{-5}$ & $7.2\times 10^{-3}$ & 269
& $7.6\times 10^{-3}$ & 1.06\\
\hline
\small$10\,000$ & 6.7$\times 10^{-6}$ & $1.0\times 10^{-3}$ & 157
& $1.1\times 10^{-3}$ & 1.04\\
\hline
\small$50\,000$ & 2.7$\times 10^{-7}$ & $9.6\times 10^{-6}$ & 36
& $9.1\times 10^{-6}$ & 0.95\\
\hline
\end{tabular}

\bigskip
The following combination of simulation and the estimate of Theorem \ref{mainthm2} seems to give efficiently good approximations for ruin probabilities. 
First fix small $\lambda_0 > 0$ and take $n_0$ such that $\lambda e^{-n_0\varphi}$ is less than $\lambda_0$. In the $j$th replication, we calculate an estimate $\hat{e}_j$ in the following way.
First apply simulation upto year $n_0$. If ruin occurs during the first $n_0$ years then put $\hat{e}_j = 1$. If ruin has not occurred then at the end of year $n_0$, the company has a random non-negative capital left. By making use of this capital, and by replacing $\lambda$ with $\lambda e^{-n_0\varphi}$,
the estimate of Theorem
\ref{mainthm2} can then be used to approximate the probability of ruin. We take $\hat{e}_j$ to be that estimate. If we have $J$ replications then the estimate of $\Prob(T<\infty)$ is the sum of the $\hat{e}_j$-observations divided by $J$. The estimation does not need much computer time since only the first $n_0$ years has to be simulated.

The above estimator consists of two parts. Firstly, the observed number of ruins divided by $J$ gives
an unbiased estimator for the probability
$\Prob(T\le n_0)$. Secondly, the sum of the estimates of Theorem \ref{mainthm2} divided by $J$ approximates the
probability $\Prob(T\in (n_0,\infty))$. This part is not unbiased but it can be expected to be accurate because $\lambda_0$ is small.

We applied the procedure by taking $\lambda = 100$ and
$\lambda_0 = 0.1$. Denote
\begin{eqnarray}
\hat{E_3} &=&  \mbox{the estimate of the ruin probability from the above procedure}.\nonumber
\end{eqnarray}
The results are given in Table 5.2.
The accuracy is measured by the quotient $\hat{E}_3/\hat{E}_2$, and it is good.

A similar ad hoc method can be used in the case where
$\xi$-variables are random. Then $n_0$ should be determined such that
$\lambda\xi_{n_0}$ is likely to be below $\lambda_0$.
Theorem \ref{mainthm2} is now applied by making use of random
$\lambda\xi_{n_0}$ instead of $\lambda$.

\section{Proofs}\label{Proofs}

\setcounter{equation}{0}
We begin by giving various lemmas to be used in the proofs of the main theorems. The proofs of the lemmas
will be given at the end of the section.

Consider first asymptotic estimates for the moments of
compound Poisson distributions.
Let ${\cal Z}, {\cal Z}_1, {\cal Z}_2,\ldots$ be an i.i.d. sequence
of non-negative random variables, and assume that
$\Prob({\cal Z}>0)>0$. Write
$${\cal S}_k = {\cal Z}_1+\cdots + {\cal Z}_{k}$$
for $k\in\nat$.
Let further ${\cal N}_{\nu}$ be a Poisson distributed random variable with the
parameter $\nu$. Assume that ${\cal N}_{\nu}$ is independent of
the ${\cal Z}$-variables, and write
\begin{equation}\label{ypLm0}
{\cal X}_{\nu} = {\cal Z}_1+\cdots + {\cal Z}_{{\cal N}_{\nu}}.
\end{equation}
Thus ${\cal X}_{\nu}$ has a compound Poisson distribution.
Let $\bar{\alpha}$ be the moment index of ${\cal Z}$, that is,
\begin{equation}\label{ypLm0a1}
\bar{\alpha} = \sup\{\alpha\ge 0\, |\, \Exp({\cal Z}^{\alpha}) < \infty\}.
\end{equation}
We will assume in the sequel that $\bar{\alpha} > 1$ so that $\Exp({\cal Z}) < \infty$.
It is well known that
\begin{equation}\label{ypLm0a2}
\limsup_{x\to\infty} (\log x)^{-1} \log\Prob({\cal Z} > x) = - \bar{\alpha}.
\end{equation}
See \shortciteN{TRHSVSJT99},
page 39. Define the function
$L_{\cal X}:(0,\bar{\alpha})\to (-\infty,\infty]$ by
\begin{equation}\label{ypLm0a3}
L_{\cal X}(\alpha) = \limsup_{\nu\to\infty}
(\log \nu)^{-1} \log \Exp\left(|{\cal X}_{\nu}-\nu\Exp({\cal Z})|^{\alpha}\right).
\end{equation}
\begin{Lm} \label{ypLm}
Assume that $\bar{\alpha}\in (1,\infty]$, and let $\alpha\in (0,\bar{\alpha})$. Then
\begin{equation}\label{pf3}
\lim_{\nu\to\infty}
(\log \nu)^{-1} \log \Exp\left({\cal X}_{\nu}^{\alpha}\right) = \alpha.
\end{equation}
Furthermore, if $0 < \alpha_1 < \alpha_2 < \bar{\alpha}$, then there exists $\varepsilon > 0$ such that
for every $\alpha\in [\alpha_1,\alpha_2]$,
\begin{equation}\label{pf3aaaa}
L_{{\cal X}}(\alpha) \le \alpha - \varepsilon.
\end{equation}
\end{Lm}

Let $\{\xi_n\}$ be a positive  process which satisfies
$(H1)-(H2)$. We next recall some basic large deviations results associated with the process.
Let $\Lambda_{\xi}^*$ be the convex conjugate of
$\Lambda_{\xi}$,
\begin{equation}\label{v14010a}
\Lambda_{\xi}^*(x) = \sup\{\alpha x -
\Lambda_{\xi}(\alpha); \alpha\in\reals\},\quad x\in\reals.
\end{equation}
Write $\theta_n = (\log\xi_n)/n$.

\begin{Lm}\label{ldLm}
Assume $(H1)$. 
Then for every $\alpha\in\reals$,
\begin{equation}\label{v1401a}
\limsup_{n\to\infty}
n^{-1}
\log \Exp\left(e^{\alpha n \theta_n}\right) =
\Lambda_{\xi}(\alpha)
\end{equation}
and for every closed set $H\subseteq\reals$,
\begin{equation}\label{v1401b}
\limsup _{n\to\infty}
n^{-1} \log\Prob\left(\theta_n\in H\right)
\le -\inf\{\Lambda_{\xi}^*(x); x\in H\}.
\end{equation}
Furthermore, if $\alpha > 0$ then
for every closed set $H\subseteq\reals$,
\begin{eqnarray}
&&\limsup _{n\to\infty}
n^{-1}
\log\Exp\left(e^{\alpha n \theta_n}{\bf 1}
(\theta_n \in H)\right)\nonumber\\
&\le& \sup\{\alpha x -\Lambda_{\xi}^*(x); x\in H\} \le \Lambda_{\xi}(\alpha).\label{v1403}
\end{eqnarray}
\end{Lm}

Consider now estimates for the distributions of the numbers of claims. 
Recall the descriptions of the $K$-variables from sections \ref{Mr} and \ref{sec21b}.
\begin{Lm}\label{cnLm}
Assume $(H1)-(H2)$. Then there exists
$\delta > 0$ such that, as $n\to\infty$,
\begin{eqnarray}
\Prob(K_n = 1) &=& \lambda \Exp(\xi_n)
+ O\left(e^{n (\Lambda_{\xi}(1)-\delta)}\right),
\label{cnLm2}\\
\Prob(K_n = 0) &=& 1-\Prob(K_n = 1)
+ O\left(e^{n (\Lambda_{\xi}(1)-\delta)}\right),\label{cnLm1}\\
\Prob(K_n \ge 2) &=& O\left(e^{n (\Lambda_{\xi}(1)-\delta)}\right)\label{cnLm3}
\end{eqnarray}
and
\begin{eqnarray}\label{cnLm4}
\Prob(K_n = 1, K_j\ge 1 \mbox{ for some }j\ge n+1)
&=& O\left(e^{n (\Lambda_{\xi}(1)-\delta)}\right).
\end{eqnarray}
Furthermore, if $\varepsilon > 0$ then there exists
$\delta > 0$ such that 
for every  $\alpha\ge 1+\varepsilon$,
\begin{eqnarray}
\Exp(K_n^{\alpha} {\bf 1}(K_n \ge 2)) &=& O\left(e^{n (\Lambda_{\xi}(1)-\delta)}\right),\quad n \to\infty.\label{cnLm3a}
\end{eqnarray}
\end{Lm}

We finally state a result associated with the process $\{{\cal Y}_{n2}\}$ of Lemma \ref{cvLm}.

\begin{Lm}\label{cvLmr}
Let $\{{\cal Y}_{n2}\}$ be a process
which satisfies \eqref{s14m10aa} and \eqref{s14m10a}
for $\alpha\in (0,\infty)$. Write
\begin{equation*}
\bar{{\cal Y}}_2 = \sup\{{\cal Y}_{n2}; n\in\nat\}\quad
\mbox{and}\quad
\underline{{\cal Y}}_2 = \inf\{{\cal Y}_{n2}; n\in\nat\}.
\end{equation*}
Then
\begin{equation}\label{pf1430}
\Exp\left(|\bar{{\cal Y}}_2|^{\alpha}\right) < \infty
\quad\mbox{and}\quad
\Exp\left(|\underline{{\cal Y}}_2|^{\alpha}\right) < \infty.
\end{equation}
\end{Lm}

{\bf Proof of Lemma \ref{cvLm}}.
Let $\alpha\in (0,\infty)$ be such that \eqref{s14m10aa} and \eqref{s14m10a} hold, and let $\bar{{\cal Y}}_2$
and $\underline{{\cal Y}}_2$ be as in Lemma \ref{cvLmr}.
By  Chebycheff's inequality,
\begin{eqnarray}\label{pf1430a}
\Prob(|\bar{{\cal Y}}_2| > u) &\le& u^{-\alpha} \Exp\left(|\bar{{\cal Y}}_2|^{\alpha}\right)
\end{eqnarray}
and
\begin{eqnarray}\label{pf1430au}
\Prob(|\underline{{\cal Y}}_2| > u) &\le& u^{-\alpha} \Exp\left(|\underline{{\cal Y}}_2|^{\alpha}\right).
\end{eqnarray}
The right hand sides of \eqref{pf1430a} and
\eqref{pf1430au} are finite by Lemma \ref{cvLmr}.
Let $\kappa, \alpha$ and $\delta\in (0,1-\kappa/\alpha)$ be such that all the conditions of
Lemma \ref{cvLm} are satisfied. Take
$\delta'$ such that
$$\delta < \delta' < 1-\kappa/\alpha,$$
and write
$$v = v(u) = u(1-u^{-\delta'}).$$
It is easy to see that
$$v(1+v^{-\delta}) \ge u(1+u^{-\delta'})$$
for large $u$ so that by \eqref{s14m3a1},

\begin{eqnarray}
\Prob\left(\bar{{\cal Y}}_1>u(1-u^{-\delta'})\right)
&=& (1+o(1)) \Prob\left(\bar{{\cal Y}}_1>
v(1+v^{-\delta})\right)\\\label{pflate}
&\le& (1+o(1))\Prob(\bar{{\cal Y}}_1> u(1+u^{-\delta'}))
= (1+o(1)) \Prob(\bar{{\cal Y}}_1>u),\quad u\to\infty.\nonumber
\end{eqnarray}
By this and \eqref{pf1430a},
\begin{eqnarray*}
\Prob(\bar{{\cal Y}}>u) &\le&
\Prob\left(\bar{{\cal Y}}_1>u(1-u^{-\delta'})\right)
+ \Prob\left(\bar{{\cal Y}}_2>u^{1-\delta'}\right)\\
&\le& (1+o(1)) \Prob(\bar{{\cal Y}}_1>u)
+ O\left(u^{-(1-\delta')\alpha}\right),\quad u\to \infty.
\end{eqnarray*}
Now $(1-\delta')\alpha > \kappa$ so that by \eqref{s14m3a0},
$\Prob(\bar{{\cal Y}}>u) \le
(1+o(1)) \Prob(\bar{{\cal Y}}_1>u)$.
On the other hand, by \eqref{s14m3a1} and \eqref{pf1430au},
\begin{eqnarray*}
\Prob(\bar{{\cal Y}}>u) &\ge&
\Prob\left(\bar{{\cal Y}}_1> u(1+u^{-\delta}),
\underline{{\cal Y}}_2 \ge -u^{1-\delta}\right)\\
&=& \Prob\left(\bar{{\cal Y}}_1>u(1+u^{-\delta})\right)
- \Prob\left(\bar{{\cal Y}}_1>u(1+u^{-\delta}),
\underline{{\cal Y}}_2 < -u^{1-\delta}\right)\\
&=& (1+o(1)) \Prob(\bar{{\cal Y}}_1>u)
+ O\left(u^{-(1-\delta)\alpha}\right).
\end{eqnarray*}
Thus $\Prob(\bar{{\cal Y}}>u) \ge
(1+o(1)) \Prob(\bar{{\cal Y}}_1>u)$.
The obtained estimates imply \eqref{s14m3a2}.
\halmos

{\bf Proof of Lemma \ref{rwlm2}.}
It is clear that $\Lambda_{\eta}'(\rate) > 0$
since $\Lambda_{\eta}(0)=0$, $\Lambda_{\eta}(\rate) > 0$
and $\Lambda_{\eta}$ is convex.
Let $b>0$ be fixed, and write
$${\cal V}'_n = ({\cal V}_n+{\cal W}){\bf 1}({\cal N}=n)
-bn{\bf 1}({\cal N}\not=n),
\quad n = 0, 1, 2,\ldots.$$
For $u>0$, write
$\tau = \tau_u = \inf\{n\in\nat\cup \{0\}\, |\, {\cal V}'_n > u\}$
($\tau = \infty$ if ${\cal V}'_n\le u$ for $n = 0, 1, 2,\ldots$). For $\varepsilon > 0$, write
$$I_u = I_{u,\varepsilon}
= [(\mu-\varepsilon)u,(\mu+\varepsilon)u].$$
Then $\{\tau = n\} = \{{\cal N}=n,
{\cal V}_n + {\cal W} > u\}$
so that
\begin{equation}\label{connections}
\{\tau < \infty\} = \{{\cal V}_{{\cal N}} + {\cal W} > u\}
\quad\mbox{ and }\quad
\{\tau\in I_u\} =
\{{\cal V}_{{\cal N}} + {\cal W}> u,\, {\cal N}\in I_u\}.
\end{equation}
Write
\begin{equation}\label{gamma}
\Gamma(\alpha) = \limsup_{n\to\infty}
n^{-1} \log \Exp\left(e^{\alpha {\cal V}'_n}\right),\quad
\alpha\in\reals.
\end{equation}
It is easy to see that for every $\alpha$
in a neighbourhood of $\rate$,
\eqref{gamma} holds as the limit and
$$\Gamma(\alpha) = \max(-\alpha b,
\Lambda_{\eta}(\alpha) -\upsilon).$$
By Theorem 2 of \shortciteN{PGWW94A}, or by theorems 3.1 and 3.2 of \shortciteN{HN94},
\begin{equation}\label{Lundberg}
\lim_{u\to\infty}
u^{-1}\log \Prob(\tau < \infty) = -\rate.
\end{equation}
Furthermore, by Theorem 4 of \shortciteN{HN95},
there exists $\varepsilon' > 0$ such that
\begin{equation}\label{typical}
\Prob(\tau\in I_u\,|\,\tau < \infty) =
1 + O(e^{-\varepsilon' u}),
\quad u\to\infty.
\end{equation}
We note that $\Gamma(\alpha)$ was assumed to be finite for  some $\alpha < 0$ in \shortciteN{HN95}, but this condition
was only needed for the sample path results of the paper.
Now \eqref{rwlm22} and \eqref{rwlm23} follow from
\eqref{connections}, \eqref{Lundberg} and \eqref{typical}.

Consider \eqref{rwlm27}. Assume first that
$$\Prob({\cal N} = n) = (1-e^{-\upsilon})e^{-\upsilon n},
\quad n = 0, 1, 2,\ldots,$$
so that ${\cal N}$ has a geometrical distribution
and $f$ is a constant function, $f(x) = 1-e^{-\upsilon}$ for $x > 0$.
Let $\zeta$ have the Bernoulli distribution with
the parameter $e^{-\upsilon}$,
$$\Prob(\zeta = 0) = 1-e^{-\upsilon},
\quad \Prob(\zeta = 1) = e^{-\upsilon}.$$
Assume that $\zeta$ is independent of everything else.
Write
$$Q = {\bf 1}(\zeta = 0) e^{{\cal W}},
\quad M = {\bf 1}(\zeta = 1) e^{\eta}
\quad\mbox{and}\quad
R = e^{\eta_1+\cdots +\eta_{{\cal N}}+  {\cal W}}.$$
Then $(Q,M)$ satisfies the conditions of Theorem \ref{basic}
with $\kappa = \rate$, and $R$ satisfies
random equation \eqref{s14m3} with this pair $(Q,M)$.
By Theorem \ref{basic},
\begin{equation}\label{c14203b}
\Prob(R > u) = (1+o(1))
\frac{\Exp(e^{\rate {\cal W}})\mu}{\rate}
(1-e^{-\upsilon}) u^{-\rate},\quad u\to\infty.
\end{equation}
This proves \eqref{rwlm27} in the case where
${\cal N}$ has a geometrical distribution.
In the general case, we make use of the
well known fact that the convergence in \eqref{rwlm21a} is uniform for $x$ in any
compact subset of $(0,\infty)$. Take $\varepsilon' > 0$ and choose
$\varepsilon > 0$ such that
$$(1-\varepsilon') f(\mu u)
\le f(xu) \le (1+\varepsilon') f(\mu u)$$
for large $u$ whenever $|x-\mu|\le\varepsilon$. Then by
\eqref{rwlm23},
\begin{eqnarray*}
\Prob({\cal V}_{{\cal N}} + {\cal W} > u)
&\le& (1+o(1)) (1+\varepsilon') f(\mu u)
\sum_{n\in I_u}e^{-\upsilon n} \Prob({\cal V}_{n} + {\cal W} > u)\\
&=&(1+o(1)) (1+\varepsilon') f(\mu u)
\sum_{n=0}^{\infty}e^{-\upsilon n}
\Prob({\cal V}_{n} + {\cal W} > u).
\end{eqnarray*}
A similar lower bound holds so that
estimate \eqref{rwlm27} for geometrically distributed
${\cal N}$ implies \eqref{rwlm27} in the general case.
\halmos

{\bf Proof of Theorem \ref{mainthm1}.}
We apply Lemma \ref{cvLm} by taking ${\cal Y}_n = Y_n$
and ${\cal Y}_{n1}$ from \eqref{pf1419}.
We begin by showing that \eqref{s14m10aa} and
\eqref{s14m10a} hold for some $\alpha > \rate_1$. Condition
\eqref{s14m10aa} does not cause any problems so that we will focus
on \eqref{s14m10a}.
Clearly,
\begin{eqnarray*}\label{pf1422c}
{\cal Y}_{n2} &=& {\cal Y}_n - {\cal Y}_{n1}\\Û
&=& \sum_{k=1}^n A_1\cdots A_{k-1} (1+i_k)
\left[V_k - \lambda m_Z \xi_k\right]\nonumber
\end{eqnarray*}
so that
\begin{equation}\label{pf1422}
{\cal Y}_{n2} - {\cal Y}_{n-1,2}
= A_1\cdots A_{n-1} (1+i_n)
\left[V_n - \lambda m_Z \xi_n\right].
\end{equation}
Choose ${\cal Z} = Z$ in Lemma \ref{ypLm}, and
take $\alpha_1 = \rate_1$ and $\alpha_2\in (\rate_1,\beta_1)$. Let $\alpha\in [\alpha_1,\alpha_2]$,
and let $H_n$ be the distribution function of $\xi_n$.
Then for any $y_0 > 0$,
\begin{eqnarray*}
\Exp\left(|V_n - \lambda m_Z \xi_n|^{\alpha}
{\bf 1}(\xi_n\ge y_0)\right)
&=&\int_{y_0}^{\infty}
\Exp\left(\left|
{\cal X}_{\lambda y}-\lambda m_Z y\right|^{\alpha}
\right) dH_n(y).\nonumber
\end{eqnarray*}
Let $\varepsilon > 0$ be
such that \eqref{pf3aaaa} holds, and take $\delta > 0$ such that
$\delta < \min(\varepsilon, \rate_1)$.
We assumed that $\Exp(\log(1+g))\ge 0$
and that $\Prob(g=0) < 1$. Hence, $\Lambda_g$ is
strictly increasing and strictly positive
on $(0,\alpha_2)$.
By Lemma \ref{ypLm}, there exist $y_0 = y_0(\alpha)$
and $c_1=c_1(\alpha)$ such that for every $n\in\nat$,
\begin{eqnarray}\label{pf1424}
\Exp\left(|V_n - \lambda m_Z \xi_n|^{\alpha}
{\bf 1}(\xi_n\ge y_0)\right)
&\le& c_1 e^{n\Lambda_g(\alpha-\delta)}.
\end{eqnarray}
We note that $y_0$ and $c_1$ depend on $\alpha$ but $\delta$ does not.
It is easy to see that
\begin{eqnarray}\label{pf1426}
\Exp\left(V_n^{\alpha} {\bf 1}(\xi_n\le y_0)\right)
\le& e^{\lambda y_0}\Exp({\cal X}_{\lambda y_0}^{\alpha}).
\end{eqnarray}
Observe that
$\Lambda_1(\rate_1) +\Lambda_g(\rate_1 -\delta)-\Lambda_g(\rate_1)< 0$
so that by continuity,
\begin{eqnarray}\label{pf1426b5}
\Lambda_1(\alpha) +\Lambda_g(\alpha -\delta)-\Lambda_g(\alpha) < 0
\end{eqnarray}
for some $\alpha\in (\rate_1,\alpha_2)$.
It follows from \eqref{pf1426} that
$\Exp\left(|V_n - \lambda m_Z \xi_n|^{\alpha}
{\bf 1}(\xi_n\le y_0)\right)$
is bounded from above by a constant.
Now $\Lambda_g(\alpha -\delta) > 0$ so that by
\eqref{pf1422} and \eqref{pf1424}, there exists a constant $c_3 = c_3(\alpha)$ such that
\begin{eqnarray}\label{pf1428}
\Exp\left(|{\cal Y}_{n2} - {\cal Y}_{n-1,2}|^{\alpha}\right)
&\le& c_3 e^{n\Lambda_A(\alpha)}
e^{n\Lambda_g(\alpha -\delta)}\nonumber\\
&=& c_3 e^{n(\Lambda_1(\alpha) +\Lambda_g(\alpha -\delta)-\Lambda_g(\alpha))}.\nonumber
\end{eqnarray}
This and \eqref{pf1426b5} imply \eqref{s14m10a}.

It is straightforward to see that under our assumptions, the conditions of Theorem
\ref{basic} are satisfied for the particular
choices of $Q$ and $M$ of \eqref{s2f15B}.
It is also clear that then $\kappa = \rate_1$ and that
$R = \bar{{\cal Y}}_{1}$
satisfies random equation \eqref{s14m3}. Apply
Theorem \ref{basic} to see that
\begin{equation}\label{pf1420}
\lim_{u\to\infty}u^{\rate_1}
\Prob(\bar{{\cal Y}}_1 > u) = C
\end{equation}
where $C$ is as in \eqref{s14m3b}.
Assume that $\Prob(q>1+s) > 0$. Then $C$ is strictly positive by \shortciteN{HN01}.
The reader is referred to Theorems 2 and 3 and to the associated discussion of the paper.
Thus all the conditions of Lemma \ref{cvLm} are satisfied
and \eqref{s2f15} holds. If $\Prob(q>1+s) = 0$ then $\bar{{\cal Y}}_{1}\le 0$ almost surely
so that $C$ of \eqref{s14m3b} equals zero. Further,
$$\Prob(T<\infty) \le \Prob(\bar{{\cal Y}}_{2} > u)$$
where $\bar{{\cal Y}}_{2}$ is as in Lemma \ref{cvLmr}.
By the same lemma and Chebycheff's inequality,
also the limit of \eqref{s2f15} equals zero.
\halmos

{\bf Proof of Proposition \ref{Kprop}.}
Estimate \eqref{cnLm09} is immediate from
\eqref{cnLm2} of Lemma \ref{cnLm} and then \eqref{cnLm00}
follows from $(H2)$. Estimate \eqref{cnLm09x} is a consequence of \eqref{cnLm00} and \eqref{cnLm3}.
Consider \eqref{cnLm11}. By Lemma \ref{ldLm},
\begin{equation}
\limsup_{n\to\infty} n^{-1}
\log \Prob(\xi_n\ge 1) \le -\inf\{\Lambda^*_{\xi}(x);
x\ge 0\}.\nonumber
\end{equation}
This proves \eqref{cnLm11} since by $(H1)$, the right hand side equals $-\infty$.
\halmos

{\bf Proof of Theorem \ref{mainthm2}.}
We will make use of Lemma \ref{cvLm} by choosing ${\cal Y}_{n} = Y_n$
and by taking ${\cal Y}_{n1}$ from
\eqref{c14203a1}.
The objective is to show that the conditions of the lemma are satisfied with $\kappa = \rate_2$. Let $p_n$, $\rho$ and $S_n$ be as described in \eqref{ss0}, \eqref{ss00}
and \eqref{ss01}. Write $W = \log ((1+i)Z)$.
Fix $\varepsilon > 0$ and let
$$J_u = \left[(\mu_2-\varepsilon)\log u,
(\mu_2+\varepsilon)\log u\right]$$
for $u > 1$. We will proceed in three steps.

{\it Step 1.} We will show that
\begin{eqnarray}
\Prob\left(\bar{{\cal Y}}_{1} > u\right)
&=& (1+o(1)) \sum_{n\in J_u}
\Prob(S_n + W > \log u)\,p_{n+1}\label{scon1}\\
&=& (1+o(1)) p \Prob(S_{\rho} + W > \log u),
\quad u\to\infty.\label{scon1x}
\end{eqnarray}
The probability
$\Prob\left(\bar{{\cal Y}}_{1} > u\right)$ can
be associated with a tail probability of a compound
distribution similarly to \eqref{lateT}. Namely, write
\begin{eqnarray*}
p'_n = \Prob(K_n = 1, K_j = 0,\,\forall j\ge n+1)\quad \mbox{and}\and\quad
p' = \sum_{n=1}^{\infty} p'_n.
\end{eqnarray*}
By Lemma \ref{cnLm}, $p'_n$ and $p_n$ are asymptotically equivalent and $p'\in (0,\infty)$.
Let $\rho'$ be a random variable such that
\begin{equation}\label{ss000}
\Prob(\rho' = n-1) = p'_n/p',
\quad n\in\nat,
\end{equation}
and assume that $\rho'$ is independent of everything else. Then
\begin{eqnarray}\label{altrep}
\Prob\left(\bar{{\cal Y}}_{1} > u\right) &=&
p'\,\Prob(S_{\rho'} + W > \log u).
\end{eqnarray}
Take ${\cal N} = \rho'$, ${\cal V}_n = S_n$ and
${\cal W} =W $, and apply Lemma \ref{rwlm2} to see that,
$$\Prob(S_{\rho'} + W> \log u)
= (1+o(1)) \Prob(S_{\rho'}+ W > \log u,\,
\rho'\in J_u),
\quad u\to\infty.$$
Hence,
$$\Prob\left(\bar{{\cal Y}}_{1} > u\right)
= (1+o(1)) \sum_{n\in J_u}
\Prob(S_n + W> \log u)\,p'_{n+1},$$
and \eqref{scon1} follows since $p_n$ and $p'_n$ are asymptotically equivalent. Take now
${\cal N} = \rho$ instead of $\rho'$, and apply Lemma \ref{rwlm2} again to see that \eqref{scon1x} holds.

{\it Step 2.}
We show that $\bar{{\cal Y}}_{1}$ satisfies conditions \eqref{s14m3a0} and \eqref{s14m3a1} of Lemma \ref{cvLm}. It follows from Lemma \ref{rwlm2}, Proposition \ref{Kprop} and \eqref{scon1x} that
\begin{equation}\label{sconxx}
\lim_{u\to\infty}
(\log u)^{-1} \log \Prob(\bar{{\cal Y}}_1>u) = -\rate_2.
\end{equation}
Thus \eqref{s14m3a0} holds with $\kappa = \rate_2$.
We will show that \eqref{s14m3a1} holds for every
$\delta > 0$. Assume first that $W\equiv 0$.
According to \eqref{scon1}, it is clear that
\begin{eqnarray}\label{scon7}
&&\Prob\left(\bar{{\cal Y}}_{1} > u(1+u^{-\delta})
\right)\\
&=& (1+o(1)) \sum_{n\in J_u}
\Prob(S_n > \log(u(1+u^{-\delta})))\,p_{n+1},\quad u\to\infty.\nonumber
\end{eqnarray}
Thus to prove \eqref{s14m3a1}, it suffices to show that
\begin{equation}\label{sconc}
\Prob(S_n > \log u)
= (1+o(1)) \Prob(S_n > \log(u(1+u^{-\delta}))),
\quad u\to\infty,
\end{equation}
uniformly for $n\in J_u$.
Let $\Lambda_{A}^*$ be the convex conjugate of
$\Lambda_{A}$,
\begin{equation}
\Lambda_{A}^*(x) = \sup\{\alpha x -
\Lambda_{A}(\alpha); \alpha\in\reals\},\quad x\in\reals.
\nonumber
\end{equation}
It follows from Theorem 1 of \shortciteN{VP65} that
for small $\varepsilon > 0$, uniformly for $n\in J_u$,
\begin{eqnarray}\label{scona}
\Prob\left(S_n > \log u\right)
&=& \Prob\left(\frac{S_{n}}{n}
> \frac{\log u}{n}\right)\\
&=& (1+o(1))
\frac{e^{-n\Lambda^*_A\left(\frac{\log u}{n}\right)}}
{\alpha_{n} \sqrt{2\pi n
\Lambda_A''(\alpha_{n})}},\quad u\to\infty,\nonumber
\end{eqnarray}
where $\alpha_n = \alpha_{n,u}$ is such that
$\Lambda_A'(\alpha_{n}) = (\log u)/n$. The probability
on the right hand side of \eqref{sconc}
is estimated similarly, and by making use of these estimates, it is easy to see by the mean value theorem
that \eqref{sconc} holds in the case where
$W\equiv 0$.

To obtain \eqref{s14m3a1} for general $W$,
it suffices by \eqref{scon1} to show that
for small $\delta > 0$,
\begin{equation}\label{gen1}
\Prob(A_1\cdots A_{\rho} (1+i)Z > u(1+u^{-\delta}))
\ge (1+o(1)) \Prob(A_1\cdots A_{\rho} (1+i)Z > u).
\end{equation}
Let $H$ be the distribution function of $(1+i) Z$.
It is easy to see that for small $\varepsilon > 0$,
\begin{equation}\label{sconxxxx}
\Prob(A_1\cdots A_{\rho} (1+i)Z > u) = (1+o(1))
\int_{u^{-\varepsilon}}^{u^{1-\varepsilon}}
\Prob\left(A_1\cdots A_{\rho} > \frac{u}{x}
\right) dH(x).
\end{equation}
To prove \eqref{gen1}, let $\delta' > 0$ be such that \eqref{s14m3a1} holds when $W\equiv 0$, and let $\varepsilon > 0$
be such that \eqref{sconxxxx} holds. Take $\delta > (1+\varepsilon)\delta'$. Then
\begin{eqnarray}\label{gen3}
&&\Prob(A_1\cdots A_{\rho} (1+i)Z > u(1+u^{-\delta}),\,
u^{-\varepsilon} \le (1+i)Z \le u^{1-\varepsilon})\\
&\ge& \int_{u^{-\varepsilon}}^{u^{1-\varepsilon}}
\Prob\left(A_1\cdots A_{\rho} > \frac{u}{x}
\left(1+\left(\frac{u}{x}\right)^{-\delta'}\right)\right)
dH(x)\nonumber\\
&=& (1+o(1)) \int_{u^{-\varepsilon}}^{u^{1-\varepsilon}}
\Prob\left(A_1\cdots A_{\rho} > \frac{u}{x}
\right) dH(x)
= (1+o(1))
\Prob(A_1\cdots A_{\rho} (1+i)Z > u).\nonumber
\end{eqnarray}
This proves \eqref{gen1}.

{\it Step 3.}
We prove that $\{{\cal Y}_{n2}\}$ satisfies conditions \eqref{s14m10aa} and \eqref{s14m10a} of Lemma \ref{cvLm}.
Condition \eqref{s14m10aa} is obviously satisfied.
Let  $\alpha\in (\rate_2,\beta_2)$. Then $\alpha > 1$. Write
\begin{eqnarray}\label{c14501a}
{\cal Y}_{n2} = {\cal Y}_n - {\cal Y}_{n1}
= {\cal W}_{n1} + {\cal W}_{n2}
\end{eqnarray}
where
\begin{eqnarray}\label{c14201a}
{\cal W}_{n1} &=&
\sum_{k=1}^n A_1\cdots A_{k-1} (1+i_k) V_k - \sum_{k=1}^n A_1\cdots A_{k-1} (1+i_k) V_k {\bf 1}(K_k = 1),
\label{c14201d}\\
{\cal W}_{n2} &=&
\sum_{k=1}^n A_1\cdots A_{k-1} (1+i_k) V_k {\bf 1}(K_k = 1)  - {\cal Y}_{n1}.\label{c14201e}
\end{eqnarray}
Let also ${\cal W}_{0j} = 0$ for $j = 1, 2$.
By Minkowski's inequality,
\begin{equation*}
\Exp\left(\left|{\cal Y}_{n2} -
{\cal Y}_{n-1,2}\right|^{\alpha}\right)
\le 2^{\alpha} \max \left\{\Exp\left(\left|{\cal W}_{nj} -
{\cal W}_{n-1,j}\right|^{\alpha}\right) ; j = 1, 2\right\}.
\end{equation*}
Thus to have \eqref{s14m10a},
it suffices to show that for some $\alpha\in (\rate_2,\beta_2)$,
\begin{equation}\label{s14m10ap}
\limsup_{n\to\infty} n^{-1}\log
\Exp\left(\left|{\cal W}_{nj} -
{\cal W}_{n-1,j}\right|^{\alpha}\right) < 0,
\quad j = 1, 2.
\end{equation}

Consider \eqref{s14m10ap} for $j=1$. Now
\begin{eqnarray}\label{c12401d1}
\Exp\left(\left|{\cal W}_{n1} -
{\cal W}_{n-1,1}\right|^{\alpha}\right)
&=& \Exp\left(\left(
A_1\cdots A_{n-1} (1+i_n) V_n {\bf 1}(K_n \ge 2)\right)^{\alpha}\right)\nonumber\\
&\le& c e^{n\Lambda_A(\alpha)}\Exp(V_n^{\alpha} {\bf 1}(K_n \ge 2))
\end{eqnarray}
where $c$ is a constant.
By Minkowski's inequality,
\begin{eqnarray}\label{c12401d2}
\Exp(V_n^{\alpha} {\bf 1}(K_n \ge 2))
&=& \sum_{h=2}^{\infty}
\Exp\left(e^{-\lambda \xi_n}\frac{(\lambda \xi_n)^h}{h!}\right)
\Exp((Z_1+\cdots +Z_h)^{\alpha})\\
&\le& \Exp(Z^{\alpha})
\sum_{h=2}^{\infty} \Exp\left(e^{-\lambda \xi_n}\frac{(\lambda \xi_n)^h}{h!}\right) h^{\alpha}
= \Exp(Z^{\alpha}) \Exp(K_n^{\alpha}{\bf 1}(K_n \ge 2)).\nonumber
\end{eqnarray}
Let $\delta > 0$ be such that \eqref{cnLm3a} holds. Then
$\Exp(V_n^{\alpha} {\bf 1}(K_n \ge 2))
= O\left(e^{n (\Lambda_{\xi}(1)-\delta)}\right)$ as
$n \to\infty$.
Take $\alpha\in (\rate_2,\beta_2)$ such that
$\Lambda_A(\alpha) + \Lambda_{\xi}(1) - \delta < 0$ to see that \eqref{s14m10ap} holds for $j=1$.

Consider \eqref{s14m10ap} for $j=2$. Now
\begin{eqnarray*}
\Exp\left(\left|{\cal W}_{n2} -
{\cal W}_{n-1,2}\right|^{\alpha}\right)
&\le& c e^{n\Lambda_A(\alpha)}
\Prob(K_n  = 1, K_j\ge 1 \mbox{ for some } j\ge n+1)
\end{eqnarray*}
where $c$ is a constant. It follows from \eqref{cnLm4}
that \eqref{s14m10ap} holds for $j=2$.

Consider now the claims of Theorem \ref{mainthm2}.
By Lemma \ref{cvLm} and steps 2 and 3, \eqref{c14203at} holds, and by step 1, $\Prob\left(\bar{{\cal Y}}_{1} > u\right)$ is asymptotically equivalent to \eqref{c14203aa}. Limit \eqref{s142030} follows from
\eqref{c14203at} and \eqref{sconxx}.
Consider \eqref{cnLm09latetyp}. Write
$T_1(u) = \inf\{n\in\nat\, ;\, {\cal Y}_{n1} > u\}$
where by convention, $T_1(u) = \infty\mbox{ if }{\cal Y}_{n1} \le u\mbox{ for every } n$.
It is seen as in \eqref{altrep} that for any $y>1$,
\begin{eqnarray*}
\Prob(T_1(u)\le y)
&=& p'\Prob(S_{\rho'}+W > \log u, \rho'\le y-1)
\end{eqnarray*}
where $p'$ and $\rho'$ are as in the first part of the proof. Thus for large $u$,
\begin{eqnarray*}
\Prob(T\le (\mu_2-\varepsilon)\log u)
&\le& \Prob\left(T_1(u/2)\le (\mu_2-\varepsilon)\log u\right)
+\Prob\left(\bar{{\cal Y}}_2 > u/2\right)\\
&\le& p'\Prob\left(S_{\rho'}+W > \log (u/2),
\rho'\le (\mu_2-\varepsilon/2)\log (u/2)\right)
+\Prob(\bar{{\cal Y}}_2 > u/2).
\end{eqnarray*}
By Lemmas \ref{rwlm2} and \ref{cvLmr},
$$\limsup_{u\to\infty} (\log u)^{-1}
\log \Prob(T\le (\mu_2-\varepsilon)\log u)
< -\rate_2.$$
For the probability
$\Prob(T\in [(\mu_2+\varepsilon)\log u,\infty))$,
the same upper bound  is obtained similarly.
Thus \eqref{cnLm09latetyp} follows from \eqref{s142030}.
Finally, \eqref{c14203aal} follows from \eqref{scon1x} and Lemma \ref{rwlm2}.
\halmos

{\bf Proof of Lemma \ref{ypLm}.} The proof of \eqref{pf3}
can be found in \shortciteN{HN10},  Lemma 3.1.
Let $\alpha\in [1,\bar{\alpha})$
and let $\varepsilon > 0$.
By Minkowski's inequality
and by \eqref{pf3},
\begin{eqnarray*}
\Exp\left(|{\cal X}_{\nu}-\nu\Exp({\cal Z})|^{\alpha}\right)
^{\frac{1}{\alpha}}
&\le& \Exp\left({\cal X}_{\nu}^{\alpha}\right)^{\frac{1}{\alpha}}
+ \nu\Exp({\cal Z})\\
&\le&(\nu^{\alpha+\varepsilon})^{\frac{1}{\alpha}}
+\nu \Exp({\cal Z})
\end{eqnarray*}
for large $\nu$.
If $\alpha\in (0,1)$ then
$(x+y)^{\alpha}\le x^{\alpha} + y^{\alpha}$
for every $x, y\ge 0$ so that
\begin{eqnarray*}
\Exp\left(|{\cal X}_{\nu}-\nu\Exp({\cal Z})|^{\alpha}\right)&\le& \nu^{\alpha+\varepsilon} + \nu^{\alpha}
\Exp({\cal Z})^{\alpha}
\end{eqnarray*}
for large $\nu$. The obtained estimates show that
$L_{\cal X}(\alpha)\le\alpha$
whenever $\alpha\in (0,\bar{\alpha})$.
By H\"older's inequality, $L_{\cal X}$ is convex.
We will show below that $L_{\cal X}(1) < 1$ so that by convexity,
$L_{\cal X}(\alpha) < \alpha$
for every $\alpha\in (0,\bar{\alpha})$.
Further, $L_{\cal X}$ is continuous so that \eqref{pf3aaaa} holds.

It remains to show that $L_{\cal X}(1) < 1$. If
$\Exp({\cal Z}^2) < \infty$ then by Schwarz's inequality,
\begin{equation}\label{pf3ce}
\Exp(|{\cal X}_{\nu} -\nu\Exp({\cal Z})|)\le\sqrt{\Var{\cal X}_{\nu}} =
\sqrt{\nu \Exp({\cal Z}^2)}.
\end{equation}
Thus $L_{\cal X}(1)\le 1/2$. In the general case, first estimate
\begin{equation}\label{pf3d}
\Exp(|{\cal X}_{\nu} -\nu\Exp({\cal Z})|)
\le\Exp(|{\cal X}_{\nu} -{\cal N}_{\nu}\Exp({\cal Z})|) +
\Exp({\cal Z}) \Exp(|{\cal N}_{\nu} -\nu|).
\end{equation}
Apply \eqref{pf3ce} with  ${\cal Z}\equiv 1$ to see that
$\Exp(|{\cal N}_{\nu} -\nu|)\le\sqrt{\nu}$.
Fix $\alpha < 2$ such that $\alpha\in (1,\bar{\alpha})$.
By H\"older's inequality and Theorem 2
of \shortciteN{BBCE65},
\begin{eqnarray}
\Exp(|{\cal Z}_1+\cdots + {\cal Z}_k -k \Exp({\cal Z})|)
&\le&
\Exp(|{\cal Z}_1+\cdots + {\cal Z}_k -k \Exp({\cal Z})|^{\alpha})
^\frac{1}{\alpha}\label{pf3f}\\
&\le& \left(2 k \Exp\left(|{\cal Z} - \Exp({\cal Z})|^{\alpha}\right)
\right)^\frac{1}{\alpha} = c k^{\frac{1}{\alpha}}\nonumber
\end{eqnarray}
where $c$ is a constant. Further,
\begin{equation}\label{pf3e}
\Exp(|{\cal X}_{\nu} -{\cal N}_{\nu}\Exp({\cal Z})|)
= \sum_{k=0}^{\infty}
e^{-\nu} \frac{\nu^k}{k!}\Exp(|{\cal Z}_1+\cdots + {\cal Z}_k -k \Exp({\cal Z})|).
\end{equation}
By \eqref{pf3f} and Jensen's inequality,
\begin{eqnarray*}
\Exp(|{\cal X}_{\nu} -{\cal N}_{\nu}\Exp({\cal Z})|)
&\le& c \Exp\left({\cal N}_{\nu}^{\frac{1}{\alpha}}\right)
\le c\nu^{\frac{1}{\alpha}}.
\end{eqnarray*}
The obtained estimates show that $L_{\cal X}(1) < 1$. \halmos

{\bf Proof of Lemma \ref{ldLm}.} The first result \eqref{v1401a} is obvious.
By $(H1)$, $\Lambda_{\xi}$ is lower semicontinuous at the origin so that \eqref{v1401b}
follows from \shortciteN{HN05a}. The first inequality of \eqref{v1403} is a special case of Varadhan's integral lemma. The proof can be found in \shortciteN{SV84} or in \shortciteN{ADOZ98}, Lemma 4.3.6,
under the additional assumption that
the level sets of $\Lambda_{\xi}^*$ are compact.
However, the proof of \shortciteN{SV84} does not need the
compactness assumption. If $\Lambda_{\xi}$ is finite in a neighbourhood of $\alpha$ then by Theorem 12.2 of \shortciteN{RR70},
\begin{eqnarray*}
\sup\{\alpha x-\Lambda_{\xi}^*(x) ; x\in\reals\} = \Lambda_{\xi}(\alpha).
\end{eqnarray*}
This proves the second inequality of \eqref{v1403}.
\halmos

{\bf Proof of Lemma \ref{cnLm}.}
We begin with some general observations.
By $(H1)$ and by convexity, $\Lambda_{\xi}$
is strictly decreasing on $(0,\infty)$.
Let $\varepsilon > 0$, $\alpha \ge 1+\varepsilon$ and
let $\delta > 0$ be small. Then
\begin{eqnarray}\label{cn1405}
\Exp\left(\xi_n^{\alpha}\right)
= O\left(e^{n (\Lambda_{\xi}(1)-\delta)}\right),
\quad n\to\infty.
\end{eqnarray}
Further, $\Prob(\xi_n > 1)\le \Exp\left(\xi_n^{\alpha}\right)$
so that
\begin{eqnarray}\label{cn1401}
\Prob(\xi_n > 1)
= O\left(e^{n (\Lambda_{\xi}(1)-\delta)}\right).
\end{eqnarray}

Consider \eqref{cnLm2}. By \eqref{cn1401},
\begin{eqnarray}
\Prob(K_n = 1) &=&
\Prob(K_n = 1, \xi_n\le 1)
+ O\left(e^{n (\Lambda_{\xi}(1)-\delta)}\right)\nonumber\\
&=&
\lambda \Exp\left(e^{-\lambda\xi_n} \xi_n {\bf 1}(\xi_n\le 1)\right) + O\left(e^{n (\Lambda_{\xi}(1)-\delta)}\right)\nonumber\\
&=& \lambda \Exp\left(\xi_n {\bf 1}(\xi_n\le 1)\right) + \lambda \Exp(\psi_n) + O\left(e^{n (\Lambda_{\xi}(1)-\delta)}\right)\nonumber
\end{eqnarray}
where
$$\psi_n = \xi_n {\bf 1}(\xi_n\le 1)\sum_{m=1}^{\infty} (-1)^m\frac{(\lambda \xi_n)^m}{m!}.$$
Thus
\begin{eqnarray}\label{ld1409}
\Prob(K_n = 1)
&=& \lambda \Exp(\xi_n)
- \lambda \Exp\left(\xi_n {\bf 1}(\xi_n > 1)\right)
+ \lambda \Exp(\psi_n) + O\left(e^{n (\Lambda_{\xi}(1)-\delta)}\right).
\end{eqnarray}
By \eqref{v1403},
\begin{eqnarray}\label{ld1411}
\limsup_{n\to\infty} n^{-1}
\log \Exp\left(\xi_n {\bf 1}(\xi_n > 1)\right)
&\le& \sup\{x-\Lambda_{\xi}^*(x) ; x\ge 0\}.
\end{eqnarray}
By $(H1)$, $\Lambda_{\xi}^*(x) = \infty$ for every $x\ge 0$ so that the supremum in \eqref{ld1411} equals
$-\infty$. Thus for any given $\delta > 0$,
\begin{eqnarray}\label{ld1414}
\Exp\left(\xi_n {\bf 1}(\xi_n > 1)\right)
&=&  O\left(e^{n (\Lambda_{\xi}(1)-\delta)}\right).
\end{eqnarray}
Clearly, $|\psi_n| \le e^{\lambda}\xi_n^2 {\bf 1}(\xi_n\le 1)$ so that by \eqref{cn1405},
$\Exp(|\psi_n|) =
O\left(e^{n (\Lambda_{\xi}(1)-\delta)}\right)$.
This together with \eqref{ld1409} and \eqref{ld1414} imply \eqref{cnLm2}.

Consider \eqref{cnLm3a}. Let
$\varepsilon > 0$ and $\alpha\ge1+\varepsilon$. By Lemma \ref{ypLm}, there exist constants $c=c(\alpha)>0$ and
$y_0 = y_0(\alpha)>0$ such that
\begin{eqnarray*}
&&\Exp(K_n^{\alpha} {\bf 1}(K_n \ge 2, \xi_n > y_0))\\
&\le& c \Exp\left(\xi_n^{2\alpha}{\bf 1}(\xi_n > y_0)\right)\le c \Exp\left(\xi_n^{1+\varepsilon}\right)
\end{eqnarray*}
for every $n\in\nat$. It follows from \eqref{cn1405} that
\begin{eqnarray}\label{c1405a}
\limsup_{n\to\infty}
n^{-1} \log \Exp(K_n^{\alpha} {\bf 1}(K_n \ge 2, \xi_n > y_0))
&\le& \Lambda_{\xi}(1) - \delta'
\end{eqnarray}
for some $\delta'>0$ which is independent of $\alpha$.
Further,
\begin{eqnarray*}
\Exp(K_n^{\alpha} {\bf 1}(K_n \ge 2, \xi_n \le y_0)) &=&
\Exp\left(e^{-\lambda\xi_n}
\sum_{m=2}^{\infty} \frac{(\lambda \xi_n)^m}{m!} m^{\alpha}{\bf 1}(\xi_n\le y_0)\right)\\
&\le&
\Exp\left(\xi_n^2\right) 
\sum_{m=2}^{\infty} \frac{\lambda^m y_0^{m-2}}{m!} m^{\alpha}
< \infty.
\end{eqnarray*}
It follows from \eqref{cn1405}
and \eqref{c1405a} that for small $\delta > 0$, \eqref{cnLm3a} holds for every $\alpha\ge 1+\varepsilon$.

Estimates \eqref{cnLm2} and \eqref{cnLm3a} imply \eqref{cnLm1} and \eqref{cnLm3}.
Thus it remains to prove \eqref{cnLm4}.
Let $\delta > 0$ be small and $c>0$ large. By \eqref{cn1401} and \eqref{cnLm3},
\begin{eqnarray*}\label{cnLm4p}
&&\Prob(K_n = 1, K_j\ge 1 \mbox{ for some }j\ge n+1)\\
&\le& \sum_{j=n+1}^{\infty} \left[\Prob(K_n = 1, K_j = 1,
\xi_n\le 1, \xi_j\le 1) +  ce^{j (\Lambda_{\xi}(1)-\delta)}\right]
+  ce^{n (\Lambda_{\xi}(1)-\delta)}\\
&=& \sum_{j=n+1}^{\infty} \Prob(K_n = 1, K_j = 1,
\xi_n\le 1, \xi_j\le 1) 
+  de^{n (\Lambda_{\xi}(1)-\delta)}
\end{eqnarray*}
where $d$ is a constant which is independent of $n$.
We conclude by Schwarz's inequality and by \eqref{cn1405} that for small $\delta > 0$,
\begin{eqnarray*}
\Prob(K_n = 1, K_j = 1, \xi_n\le 1, \xi_j\le 1)
&=&\Exp\left(e^{-\lambda\xi_n} \lambda\xi_n\, e^{-\lambda\xi_j} \lambda\xi_j
{\bf 1}(\xi_n\le 1){\bf 1}(\xi_j\le 1)\right)\\
&\le& \lambda^2 \Exp\left(\xi_n^2 {\bf 1}(\xi_n\le 1)\right)
^\frac{1}{2}
\Exp\left(\xi_j^2 {\bf 1}(\xi_j\le 1)\right)
^\frac{1}{2}\\
&\le& \lambda^2
e^{n (\Lambda_{\xi}(1)-\delta)}
e^{\frac{1}{2}(j-n) (\Lambda_{\xi}(1)-\delta)}
\end{eqnarray*}
for every $j\ge n+1$ for large $n$.
The obtained estimates imply \eqref{cnLm4}.
\halmos

{\bf Proof of Lemma \ref{cvLmr}}.
We only prove the first
inequality of \eqref{pf1430}.
Obviously,
\begin{equation}
|\bar{{\cal Y}}_2| \le \sum_{n=1}^{\infty}
|{\cal Y}_{n2} - {\cal Y}_{n-1,2}|.
\end{equation}
Suppose that $\alpha\in (0,1)$. Then
$(x+y)^{\alpha}\le x^{\alpha}+y^{\alpha}$
for every $x, y\ge 0$. Thus
\begin{eqnarray*}
\Exp\left(|\bar{{\cal Y}}_2|^{\alpha}\right)
&\le& \sum_{n=1}^{\infty}
\Exp\left(|{\cal Y}_{n2} - {\cal Y}_{n-1,2}|^{\alpha}\right).
\end{eqnarray*}
The terms of the series are finite by \eqref{s14m10aa}, and
by \eqref{s14m10a}, there exists $\delta > 0$ such that
\begin{equation}\label{s14m10aaa}
\Exp\left(|{\cal Y}_{n2} - {\cal Y}_{n-1,2}|^{\alpha}\right)
\le e^{-n\delta}
\end{equation}
for large $n$. Thus 
$\Exp\left(|\bar{{\cal Y}}_2|^{\alpha}\right) <\infty$.
Let now $\alpha\ge 1$. By Minkowski's inequality,
\begin{eqnarray*}
\Exp\left(|\bar{{\cal Y}}_2|^{\alpha}\right)^{\frac{1}{\alpha}}
&\le& \sum_{n=1}^{\infty}
\Exp\left(|{\cal Y}_{n2} - {\cal Y}_{n-1,2}|^{\alpha}\right)
^{\frac{1}{\alpha}}.
\end{eqnarray*}
The right hand side is finite by \eqref{s14m10aa} and
\eqref{s14m10aaa} so that
$\Exp\left(|\bar{{\cal Y}}_2|^{\alpha}\right) <\infty$.
\halmos

\renewcommand{\baselinestretch}{1.0}
{\small

\bibliographystyle{chicago}
\bibliography{Bbl}

\begin{thebibliography}{}

\bibitem[\protect\citeauthoryear{Asmussen and Kl\"{u}ppelberg}{Asmussen and
  Kl\"{u}ppelberg}{1996}]{SACK96}
Asmussen, S. and C.~Kl\"{u}ppelberg (1996).
\newblock Large deviations results for subexponential tails, with applications
  to insurance risk.
\newblock {\em Stoch. Proc. Appl.\/}~{\em \bf 64}, 103--125.

\bibitem[\protect\citeauthoryear{Daykin, Pentik\"ainen, and Pesonen}{Daykin
  et~al.}{1994}]{CDTPMP94}
Daykin, C.~D., T.~Pentik\"ainen, and M.~Pesonen (1994).
\newblock {\em Practical Risk Theory for Actuaries}.
\newblock London: Chapman \& Hall.

\bibitem[\protect\citeauthoryear{Dembo and Zeitouni}{Dembo and
  Zeitouni}{1998}]{ADOZ98}
Dembo, A. and O.~Zeitouni (1998).
\newblock {\em Large Deviations Techniques and Applications\/} (2nd ed.).
\newblock Berlin: Springer--Verlag.

\bibitem[\protect\citeauthoryear{Embrechts, Maejima, and Teugels}{Embrechts
  et~al.}{1985}]{PEMMJT85}
Embrechts, P., M.~Maejima, and J.~Teugels (1985).
\newblock Asymptotic behaviour of compound distributions.
\newblock {\em Astin Bulletin\/}~{\em \bf 15}, 45--48.

\bibitem[\protect\citeauthoryear{Glynn and Whitt}{Glynn and
  Whitt}{1994}]{PGWW94A}
Glynn, P.~W. and W.~Whitt (1994).
\newblock Logarithmic asymptotics for steady-state probabilities in a
  single-server queue.
\newblock {\em J. Appl. Prob.\/}~{\em \bf 31A}, 131--156.

\bibitem[\protect\citeauthoryear{Goldie}{Goldie}{1991}]{CG91}
Goldie, C.~M. (1991).
\newblock Implicit renewal theory and tails of solutions of random equations.
\newblock {\em Ann. Appl. Probab.\/}~{\em \bf 1}, 126--166.

\bibitem[\protect\citeauthoryear{Grandell}{Grandell}{1997}]{JG97}
Grandell, J. (1997).
\newblock {\em Mixed Poisson Processes}.
\newblock Chapman \& Hall, London.

\bibitem[\protect\citeauthoryear{Norberg}{Norberg}{1993}]{RN93}
Norberg, R. (1993).
\newblock Prediction of outstanding liabilities in non-life insurance.
\newblock {\em Astin Bulletin\/}~{\em \bf 23}, 95--115.

\bibitem[\protect\citeauthoryear{Nyrhinen}{Nyrhinen}{1994}]{HN94}
Nyrhinen, H. (1994).
\newblock Rough limit results for level crossing probabilities.
\newblock {\em J. Appl. Probab.\/}~{\em \bf 31}, 373--382.

\bibitem[\protect\citeauthoryear{Nyrhinen}{Nyrhinen}{1995}]{HN95}
Nyrhinen, H. (1995).
\newblock On the typical level crossing time and path.
\newblock {\em Stoch. Proc. Appl.\/}~{\em \bf 58}, 121--137.

\bibitem[\protect\citeauthoryear{Nyrhinen}{Nyrhinen}{2001}]{HN01}
Nyrhinen, H. (2001).
\newblock Finite and infinite time ruin probabilities in a stochastic economic
  environment.
\newblock {\em Stoch. Proc. Appl.\/}~{\em \bf 92}, 265--285.

\bibitem[\protect\citeauthoryear{Nyrhinen}{Nyrhinen}{2005}]{HN05a}
Nyrhinen, H. (2005).
\newblock Upper bounds of the {G}\"artner-{E}llis theorem for the sequences of
  random variables.
\newblock {\em Statist. Probab. Lett.\/}~{\em \bf 73}, 57--60.

\bibitem[\protect\citeauthoryear{Nyrhinen}{Nyrhinen}{2010}]{HN10}
Nyrhinen, H. (2010).
\newblock Economic factors and solvency.
\newblock {\em Astin Bulletin\/}~{\em \bf 40}, 889--916.

\bibitem[\protect\citeauthoryear{Paulsen}{Paulsen}{2008}]{JP08}
Paulsen, J. (2008).
\newblock Ruin models with investment income.
\newblock {\em Probab. Surv.\/}~{\em \bf 5}, 416--434.

\bibitem[\protect\citeauthoryear{Pentik\"ainen, Bonsdorff, Pesonen, Rantala,
  and Ruohonen}{Pentik\"ainen et~al.}{1989}]{TPHBMPJRMR89}
Pentik\"ainen, T., H.~Bonsdorff, M.~Pesonen, J.~Rantala, and M.~Ruohonen
  (1989).
\newblock {\em Insurance Solvency and Financial Strength}.
\newblock Finnish Insurance Training and Publishing Company, Helsinki.

\bibitem[\protect\citeauthoryear{Pentik\"ainen and Rantala}{Pentik\"ainen and
  Rantala}{1982}]{TPJR82}
Pentik\"ainen, T. and J.~Rantala (1982).
\newblock {\em Solvency of Insurers and Equalization Reserves, Vol. I and II}.
\newblock The Insurance Publishing Company, Helsinki.

\bibitem[\protect\citeauthoryear{Petrov}{Petrov}{1965}]{VP65}
Petrov, V.~V. (1965).
\newblock On the probabilities of large deviations for sums of independent
  random variables.
\newblock {\em Theory Probab. Appl.\/}~{\em \bf 10}, 287--298.

\bibitem[\protect\citeauthoryear{Rantala}{Rantala}{1984}]{JR84}
Rantala, J. (1984).
\newblock {\em An application of stochastic control theory to insurance
  business}.
\newblock PhD. Thesis. University of Tampere.

\bibitem[\protect\citeauthoryear{Rockafellar}{Rockafellar}{1970}]{RR70}
Rockafellar, R.~T. (1970).
\newblock {\em Convex Analysis}.
\newblock Princeton: Princeton Univ. Press.

\bibitem[\protect\citeauthoryear{Rolski, Schmidli, Schmidt, and Teugels}{Rolski
  et~al.}{1999}]{TRHSVSJT99}
Rolski, T., H.~Schmidli, V.~Schmidt, and J.~Teugels (1999).
\newblock {\em Stochastic Processes for Insurnace and Finance}.
\newblock Chichester, UK: Wiley.

\bibitem[\protect\citeauthoryear{Ruohonen}{Ruohonen}{1988}]{ICA88}
Ruohonen, M. (1988).
\newblock The claims occurrence process and the {I.B.N.R}. problem.
\newblock In {\em Proceedings of international congress of actuaries, Part 4},
  Helsinki, pp.\  113--123. Gummerus Oy, Jyv\"askyl\"a.

\bibitem[\protect\citeauthoryear{Teugels}{Teugels}{1985}]{JT85}
Teugels, J. (1985).
\newblock Approximation and estimation of some compound distributions.
\newblock {\em Insurance: Mathematics and Economics 4\/}, 143--153.

\bibitem[\protect\citeauthoryear{Varadhan}{Varadhan}{1984}]{SV84}
Varadhan, S. R.~S. (1984).
\newblock {\em Large Deviations and Applications}.
\newblock Philadelphia: SIAM.

\bibitem[\protect\citeauthoryear{von Bahr and Esseen}{von Bahr and
  Esseen}{1965}]{BBCE65}
von Bahr, B. and C.~Esseen (1965).
\newblock Inequalities for the $r$th absolute moment of a sum of random
  variables, $1\le r\le 2$.
\newblock {\em Ann. Math. Statist.\/}~{\em \bf 36}, 299--303.

\end{thebibliography}

\vspace{.5cm}
\noindent
{\sc Harri Nyrhinen\\
Department of Mathematics and Statistics\\
P.O. Box 68 (Gustaf H\"{a}llstr\"{o}min Katu 2b)\\
FIN 00014, University of Helsinki, Finland}\\
{\it E-mail:}  harri.nyrhinen@helsinki.fi
}

\end{document}